\documentclass[11pt,reqno]{amsart}
\pdfoutput=1
\usepackage[utf8]{inputenc}

\usepackage[dvips]{graphicx}

\usepackage{amssymb, amsmath, amsthm, pstricks}
\usepackage{bbm} 

\makeatletter
\@namedef{subjclassname@2020}{\textup{2020} Mathematics Subject Classification}
\makeatother

\usepackage[hyphens]{url}

\usepackage[backend=biber,bibencoding=utf8,style=alphabetic
,sorting=nyt,firstinits=true,doi=true,isbn=false,backref=true,url=true,hyperref=true,maxbibnames=99,maxcitenames=5]{biblatex}
\renewbibmacro{in:}{}
\DeclareFieldFormat{pages}{#1}
\usepackage{hyperref}

\usepackage{enumitem}

\newcommand{\cf}{{\it cf.}}
\newcommand{\ie}{{\it i.e.}}
\newcommand{\eg}{{\it e.g.}}

\newcommand{\term}{\emph} 

\DeclareMathOperator\Img{Im} 
\DeclareMathOperator\Int{Int} 
\newcommand{\diff}{\mathop{}\!\mathrm d} 
\DeclareMathOperator\length{length}
\DeclareMathOperator\area{area}
\DeclareMathOperator\vol{vol}

\DeclareMathOperator\inj{inj} 
\DeclareMathOperator\conv{conv} 
\DeclareMathOperator\Prob{\mathbb P}
\DeclareMathOperator\dProb{\diff\mathbb P}
\DeclareMathOperator\Expe{\mathbb E}
\DeclareMathOperator\Leg{\mathcal L} 
\newcommand{\convol}{*} 
\newcommand{\N}{{\mathbb N}}
\newcommand{\Z}{{\mathbb Z}}

\newcommand{\R}{{\mathbb R}}
\newcommand{\C}{\mathcal{C}}
\newcommand{\e}{\mathrm{e}}
\newcommand{\ii}{\mathrm{i}}
\newcommand{\eps}{\varepsilon} 
\newcommand{\II}{{\mathcal I}} 
\newcommand{\DD}{\mathcal{D}} 
\newcommand{\W}{\mathcal{W}} 
\newcommand{\HH}{\mathcal{H}} 
\newcommand{\strip}{\mathcal{S}}
\newcommand{\ext}{\textrm{ext}} 

\long\def\forget#1\forgotten{}

\numberwithin{equation}{section}
\newtheorem{theorem}{Theorem}[section]
\newtheorem{proposition}[theorem]{Proposition}
\newtheorem{corollary}[theorem]{Corollary}
\newtheorem{lemma}[theorem]{Lemma}
\newtheorem{claim}[theorem]{Claim}
\theoremstyle{definition}
\newtheorem{definition}[theorem]{Definition}
\newtheorem{example}[theorem]{Example}
\newtheorem{remark}[theorem]{Remark}

\begin{document}

\title
{Minimal~area~of~Finsler~disks with~minimizing~geodesics}

\author[M.~Cossarini]{Marcos Cossarini}
\author[S.~Sabourau]{St\'ephane Sabourau}

\thanks{Partially supported by the B\'ezout Labex (ANR-10-LABX-58) and the ANR project Min-Max (ANR-19-CE40-0014).}

\address{\parbox{\linewidth}{Univ Paris Est Creteil, CNRS, LAMA, F-94010 Creteil, France \\
Univ Gustave Eiffel, LAMA, F-77447 Marne-la-Vall\'ee, France}}

\email{marcos.cossarini@u-pec.fr}

\email{stephane.sabourau@u-pec.fr}



\begin{abstract}
We show that the Holmes--Thompson area of every Finsler disk of radius~$r$ whose interior geodesics are length-minimizing is at least~$\frac{6}{\pi} r^2$.
Furthermore, we construct examples showing that the inequality is sharp and observe that the equality case is attained by a non-rotationally symmetric metric.
This contrasts with Berger's conjecture in the Riemannian case, which asserts that the round hemisphere is extremal.
To prove our theorem we discretize the Finsler metric using random geodesics. As an auxiliary result, we include a proof of the integral geometry formulas of Blaschke and Santal\'o for Finsler manifolds with almost no trapped geodesics.
\end{abstract}

\maketitle

\tableofcontents

\section{Introduction}

Isoembolic inequalities on Riemannian manifolds are curvature-free volume estimates in terms of the injectivity radius.
The first sharp isoembolic inequality valid in all dimension is due to Berger~\cite{berger1980borne} who showed that the volume of every closed Riemannian $n$-manifold~$M$ satisfies
\begin{equation} \label{eq:berger}
\vol(M) \geq \alpha_n \left( \frac{\inj(M)}{\pi} \right)^n
\end{equation}
where~$\alpha_n$ is the volume of the canonical $n$-sphere.
Furthermore, equality holds if and only if~$M$ is isometric to a round sphere.
The two-dimensional case was proved earlier in~\cite{berger1976relations}.

A long standing conjecture in Riemannian geometry also due to Berger asserts that every ball~$B(r)$ of radius $r \leq \frac{1}{2} \inj(M)$ in a closed Riemannian $n$-manifold~$M$ satisfies
\begin{equation}\label{eq:conj-berger}
\vol B(r) \geq \frac{\alpha_n}{2} \left( \frac{2r}{\pi} \right)^n
\end{equation}
with equality if and only if~$B(r)$ is isometric to a round hemisphere of (intrinsic) radius~$r$.
This can be viewed as a local version of the sharp isoembolic inequality~\eqref{eq:berger}.
This conjecture is open even in the two-dimensional case where the previous inequality can be written as
\begin{equation} \label{eq:conj2}
\area \, D(r) \geq \frac{8}{\pi} \, r^2.
\end{equation}

An account on isoembolic inequalities and Berger's conjecture is given in~\cite[\S6]{croke2003universal}.
A non-sharp volume estimate~$\vol \, B(r) \geq c_n r^n$ was established by Berger \cite{berger1976relations}, \cite{berger1977volume} for~$n=2$ or~$3$, and by Croke~\cite[Proposition~14]{croke1980isoperimetric} for every~$n$.
The conjecture (with a sharp constant) is satisfied for metrics of the form
$\diff s^2 = \diff r^2 +f(r,\theta)^2 \, \diff \theta^2$ in polar coordinates when~$n \geq 3$; see~\cite{croke1983volume}.
In~\cite{croke1984curvature}, Croke also showed that the optimal inequality~\eqref{eq:conj2} holds true on average over all balls~$B(r)$ of~$M$.
In the two-dimensional case, the best general estimate
$\area \, D(r) \geq \frac{8-\pi}{2} r^n$ can be found in~\cite{croke2009area}.
The lower bound~\eqref{eq:conj2} on the area of~$D(r)$ has recently been obtained in~\cite{chambers2017area} by Chambers--Croke--Liokumovich--Wen under the stronger hypothesis that~$r \leq \frac{1}{2} \conv(M)$, where~$\conv(M)$ is the convexity radius of~$M$.
(This implies that ~$r \leq \frac{1}{4} \inj(M)$, since $\conv(M) \leq \frac{1}{2} \inj(M)$.)
Note, however, that this stronger condition rules out the possibility that $D(r)$ is a hemisphere of intrinsinc radius~$r$, which is the only expected equality case of~\eqref{eq:conj2}.

The condition that $r \leq \frac{1}{2} \inj(M)$ in Berger's conjecture~\eqref{eq:conj-berger} can be relaxed by requiring instead that every interior geodesic in~$B(r)$ is length-minimizing.
The results of \cite{berger1976relations}, \cite{berger1977volume} and \cite[Proposition~14]{croke1980isoperimetric}, for instance, still hold under this more general condition.

\medskip

In this article, we consider the case of disks with a self-reverse Finsler metric whose interior geodesics are length-minimizing.
(A precise definition of Finsler metrics and area can be found in Section~\ref{sec:finsler}.)
It is natural to expect that the inequality~\eqref{eq:conj2} holds in this setting.
This is the case for isosystolic inequalities on the projective plane, where the canonical round metric minimizes the systolic area among both Riemannian and Finsler metrics; see~\cite{ivanov2002two}, \cite{ivanov2011filling}.
However, we show that the round hemisphere is not area minimizing among Finsler metric disks of the same radius whose interior geodesics are length-minimizing.
More precisely, we establish a sharp isoembolic inequality for Finsler metrics in the two-dimensional case under the assumption that every interior geodesic is length-minimizing.
We observe that the extremal metric is not Riemannian and, surprisingly, not even rotationally symmetric.

\medskip

Before stating our main result, let us introduce the following definition.

\begin{definition}
A \term{Finsler disk~$D$ of radius~$r$ with minimizing interior geodesics} is a disk with a Finsler metric such that
\begin{itemize}
\item every interior point of~$D$ is at distance less than~$r$ from a specified center point~$O$;
\item every point of~$\partial D$ is at distance exactly~$r$ from~$O$;
\item every interior geodesic of~$D$ is length-minimizing.
\end{itemize}
For instance, a ball of radius~$r$ on a complete Finsler plane with no conjugate points is a Finsler disk of radius~$r$ with minimizing interior geodesics.
\end{definition}

The optimal version of Berger's conjecture for Finsler surfaces with self-reverse metric is given by the following result.
We emphasize that we make no assumptions on the convexity radius.

\begin{theorem} \label{theo:6pi-intro}
Let~$D$ be a self-reverse Finsler metric disk~$D$ of radius~$r$ with minimizing interior geodesics.
Then the Holmes--Thompson area of~$D$ satisfies
\[  \area(D) \geq \frac{6}{\pi} \, r^2.  \]
Furthermore, the inequality is optimal.
\end{theorem}

The lower bound is attained by a non-smooth space consisting of a disk of radius~$r$ centered at the tip of the cone obtained by gluing together three copies of a quadrant of the $\ell^1$-plane.
(Recall that the $\ell^1$-plane is the normed plane where unit balls have the least possible area, according to Mahler's theorem on convex bodies in the plane.) 
Note that this disk is not rotationally symmetric.
In Section~\ref{sec:extremal2} we use Busemann's construction of projective metrics (developed in relation with Hilbert's fourth problem) to give another description of this non-smooth extremal metric.
More precisely, we define a non-smooth projective metric on the plane where the disk of radius $r$ centered at the origin has area $\frac 6\pi r^2$.
Then, we approximate this non-smooth projective metric by smooth projective metrics (which are therefore Finsler and have minimizing interior geodesics) where the area of the disk converges to $\frac 6\pi r^2$, proving that the inequality of Theorem~\ref{theo:6pi-intro} is sharp.

\medskip

Let us further comment on the result proved in~\cite{chambers2017area} for Riemannian disks~$D(r) \subseteq M$ of radius $r \leq \frac{1}{2} \conv(M)$.
As previously mentioned, this excludes the possibility that $D(r)$ is a hemisphere of intrinsinc radius~$r$.
Still, the argument in~\cite{chambers2017area} is valid for Finsler surfaces with self-reverse metric, except for the proof of their Lemma~2.1, which is purely Riemannian.
Therefore, a self-reverse Finsler metric disk of radius~$r$ in which the distance function from each given point is convex along all geodesics satisfies~\eqref{eq:conj2}. 
The extremal surfaces that we construct in this paper violate this inequality, however this poses no contradiction because they have a vanishing convexity radius.

\medskip

Instead of the Holmes--Thompson area, one could consider the Busemann--Hausdorff area, which, in general, is bounded below by the former; see~\cite{duran1998volume}.
However, the Busemann--Hausdorff area of the extremal metric in Theorem~\ref{theo:6pi-intro} is equal to~$\frac{3}{4} \pi r^2$, which is greater than the area of the round hemisphere of intrinsic radius~$r$ that is conjectured to be minimal.

\medskip

The proof of Theorem~\ref{theo:6pi-intro} and the construction of extremal and almost extremal metrics occupy the whole article.
The approach, based on a discretization of the metric (\cf~\cite{cossarini2018discrete}), is fairly robust and new in this context.

\medskip

The article is organized as follows.

In Section~\ref{sec:finsler}, we recall the notions of Finsler manifolds, their Holmes--Thompson measure, and their geodesics described from the Hamiltonian point of view.

In Section~\ref{sec:integral}, we go over the standard proofs of the integral geometry formulas of Blaschke and Santaló, showing that they are valid for Finsler manifolds with almost no trapped geodesics.
In the case of a disk as in Theorem~\ref{theo:6pi-intro}, the formulas say that the length of a curve in the disk is proportional to the expected number of intersections with a random geodesic, and the area of a region is proportional to the expected length of the intersection with a random geodesic.

In Section~\ref{sec:discretization}, we introduce the notion of a quasi wall system on a surface, generalizing the wall systems studied in~\cite{cossarini2018discrete}.
A quasi wall system on a surface is a 1-dimensional submanifold satisfying certain conditions.
It determines a discrete metric, according to which
the length of a curve is its number of intersections with the quasi wall system,
and the area of the surface is the number of self-intersections of the quasi wall system.
We show how to approximate a self-reverse Finsler metric with minimizing geodesics by a quasi wall system consisting of random geodesics.
To prove the approximation properties we use the integral geometry formulas to compute the expected values of discrete length and area, and then we apply the law of large numbers.

In Section~\ref{sec:discrete-to-smooth}, we use this approximation result to show that Theorem~\ref{theo:6pi-intro} follows from an analogous theorem on simple discrete metric disks.

Sections~\ref{sec:semi-circle}, \ref{sec:inadmissible}, \ref{sec:adjacent} and~\ref{sec:proof} are devoted to the proof of this discrete theorem. 
The proof is based on identifying certain configurations on a quasi wall system and operating on these configurations in order to transform a simple discrete disk into a new one of less area.
When the disk has minimum area, none of these configurations is present, and this implies that the quasi wall system is of a special kind where we can compute a lower bound for the area.

In Section~\ref{sec:extremal1}, we construct a simple discrete disk of minimal discrete area and show that it is unique up to isotopy.

In Section~\ref{sec:extremal2}, we use Busemann's construction of projective metrics to obtain continuous versions of our discrete area-minimizing disk.

Finally, Section~\ref{sec:regularity} is an appendix where we show that on a Finsler surface with boundary, distance-realizing curves are $C^1$.
\\

\noindent {\it Acknowledgment.}
The first author thanks the Laboratoire d'analyse et de math\'ematiques appliqu\'ees at the Universit\'e Gustave Eiffel/Universit\'e Paris Est Cr\'eteil and the Groupe Troyanov at the \'Ecole Polythechnique F\'ed\'erale de Lausanne for hosting him as a postdoc while this work was done.
The second author would like to thank the Fields Institute and the Department of Mathematics at the University of Toronto, where part of this work was accomplished, for their hospitality.
The authors thank the referees for their comments, which helped improve the exposition.

\section{Finsler metrics and Holmes--Thompson volume} \label{sec:finsler}

In this section, we recall basic definitions of Finsler geometry.

\subsection{Finsler metrics}

Let us recall the definition of a Finsler metric. 

\begin{definition}
A \term{Finsler metric} on a smooth manifold~$M$ is a continuous function $F:TM \to [0,\infty)$ on the tangent bundle~$TM$ of~$M$ satisfying the following properties (here, $F_x:=F|_{T_x M}$ for short):
\begin{enumerate}
\item\label{def:positive_homogeneity}
Positive homogeneity: $F_x(tv) = t \, F_x(v)$ for every $v \in T_x M$ and $t \geq 0$.
\item Subadditivity: $F_x(v+w)\leq F_x(v) + F_x(w)$ for every $v,w\in T_xM$.
\item Positive definiteness: $F_x(v)>0$ for every nonzero $v \in T_xM$. \label{def:positiveness}
\item\label{def:homogeneous_smoothness}
Smoothness: $F$ is smooth outside the zero section.
\item\label{def:strong_convexity}
Strong convexity: for any two linearly independent vectors $v,w\in T_xM$, the Hessian value $q_v(w)=\left.\frac{\diff^2}{\diff t^2}\right|_{t=0}F_x(v+tw)$ is strictly positive.
\end{enumerate}
Additionally, a Finsler metric~$F$ may be or not be
\begin{enumerate}[resume]
\item Self-reverse: $F_x(-v)=F_x(v)$ for every $v\in T_xM$.
\end{enumerate}
\end{definition}

Equivalently, one could define a Finsler metric by replacing~\eqref{def:positiveness} and~\eqref{def:strong_convexity} with the condition that for every nonzero vector~$v \in TM$, the Hessian of~$F^2$ at~$v$ is positive definite; see~\cite{cossarini2020finsler}.


In each tangent space~$T_xM$, the unit ball and unit sphere determined by the norm~$F_x$ are \[
B_xM = \{ v\in T_xM \mid F_x(v) \leq 1 \} \text{ and }
U_xM = \{ v\in T_xM \mid F_x(v) =    1 \}.  \]
Similarly, in the cotangent space~$T_x^*M$, the norm~$F_x^*$ dual to~$F_x$ determines a unit co-ball~$B_x^*M$ and a unit co-sphere~$U_x^*M$.

\begin{remark}\label{rmk:extend}
To handle technical details in case~$M$ has nonempty boundary, we extend the metric~$F$ to a manifold $M^+\supseteq M$, of the same dimension as~$M$ but without boundary.
\end{remark}

\subsection{Length, geodesics and distance-realizing arcs} \label{subsec:length}
\begin{definition} Let~$M$ be a manifold with a Finsler metric~$F$.
The \term{length} of a piecewise-$C^1$ curve $\gamma:I\to M$ is defined as the integral of its \term{speed}~$F(\gamma'(t))$, that is,
\begin{equation} \label{eq:length}
\length(\gamma) = \int_I F(\gamma'(t)) \diff t
\end{equation}
and the \term{distance} $d_F(x,y)$ between two points~$x$ and~$y$ in~$M$ is the infimum length of a curve~$\gamma$ in~$M$ joining~$x$ to~$y$.

A \term{distance-realizing curve} is a curve $\gamma:I\to M$ such that
\[ d_F(\gamma(t),\gamma(t')) = t' - t \]
for every~$t<t'$.

A \term{geodesic} of~$M$ is a smooth, unit-speed curve $\gamma:I\to M$ that is extremal for the length functional.
In case~$M$ has boundary, the extremality is defined by considering variations in~$M^+$, see Remark~\ref{rmk:extend}.
Thus the geodesics of~$M$ are the geodesics of~$M^+$ that are contained in~$M$.
Equivalently, the geodesics of~$M$ are the unit-speed curve curves that satisfy the Euler--Lagrange equation for the length functional; see Definition~\ref{def:cogeodesic} below for an explicit equation in terms of momentum.
\end{definition}

In a compact connected Finsler manifold, every pair of points are joined by a distance-realizing arc.
\footnote{A proof for more general, complete self-reverse metrics is given~\cite[\S1.12]{gromov2007metric}; see also~\cite[Theorem~9.1]{mennucci2014geodesics} for directed metrics.} 
A distance-realizing arc contained in the interior of~$M$ is necessarily a geodesic and is therefore smooth. However, a distance-realizing arc of~$M$ does not necessarily lie in the interior of~$M$, even if its endpoints do.
Still, if the manifold is two-dimensional, then every distance-realizing arc is~$C^1$ and has unit speed; see Theorem~\ref{thm:shortest_are_C1}.
Thus, in a compact Finsler surface, every pair of points~$x,\,y$ are joined by a~$C^1$ arc of length~$d_F(x,y)$.


\subsection{Symplectic structure on the cotangent bundle}

Recall also some definitions about the geodesic flow of a Finsler manifold from the Hamiltonian viewpoint; see~\cite[Chap. 7--9]{arnold1989mathematical}, 
\cite{paiva2006problems} and \cite{cossarini2020finsler}.

\begin{definition} \label{def:tautological}
Let~$M$ be a manifold.
The \term{tautological one-form}~$\alpha_M$ on~$T^*M$ is defined as
\[  \alpha_M|_\xi(V) = \xi (\diff \pi_\xi(V))  \]
for every $\xi \in T^*M$ and $V \in T_\xi T^*M$, where $\pi:T^*M \to M$ is the canonical projection.
The \term{standard symplectic form}~$\omega_M$ on~$T^*M$ is given by
\[  \omega_M = \diff \alpha_M.  \]
Using canonical coordinates $(x_i,\xi_i)$ on~$T^*M$, 
these forms can be expressed as
\begin{equation}\label{eq:alpha_omega_canon}
\alpha_M = \sum_i \xi_i \diff x_i,\qquad
\omega_M = \sum_i \diff \xi_i \wedge \diff x_i.
\end{equation}
\end{definition}

\begin{definition}\label{def:legendre}
Let~$(M,F)$ be a Finsler manifold. The \term{Legendre map} 
\[  \Leg: UM \to U^*M,  \]
is defined as follows: the image of a unit vector $v \in U_xM$ is the unique unit covector~$\xi \in U_x^*M$ such that~$\xi(v)=1$.
Since~$F$ is strongly convex, the Legendre map is a diffeomorphism. Its inverse is the Legendre map associated to the dual metric~$F^*$ on~$T^*M$, which is also strongly convex.
The unit covectors will also be referred to as \term{momentums}.
The \term{Hamiltonian lift} of a unit-speed curve~$\gamma$ in~$M$ is the curve $t\mapsto \Leg(\gamma'(t))$ in~$U^*M$.
\end{definition}

\begin{definition}\label{def:cogeodesic}
The \term{cogeodesic vector field} of a Finsler manifold $M$ is the vector field~$Z$ on~$U^*M$ given by the equations
\begin{align*}
\iota_Z({\omega_M} |_{U^*M}) &= 0\\
\iota_Z(\alpha_M)           &= 1.
\end{align*}
where $\iota_Z$ is the operator that contracts a differential form with the vector field $Z$.
The integral curves of $Z$ are the Hamiltonian lifts of the geodesics in~$M$; see~\cite{cossarini2020finsler}.
\end{definition}

It follows from the Cartan formula that the forms ~$\alpha_M$ and~$\omega_M$ restricted to~$U^*M$ are invariant under the cogeodesic flow.

\subsection{Holmes--Thompson volume}
We will consider the following notion of volume.

\begin{definition}
The \term{Holmes--Thompson volume} of a Finsler $n$-manifold~$M$ is defined as the symplectic volume of its unit co-ball bundle~$B^*M\subseteq T^*M$, divided by the volume~$\epsilon_n$ of the Euclidean unit ball in~$\R^n$. That is,
\begin{equation} \label{eq:volHT}
\vol(M)=\frac 1{\epsilon_n}\int_{B^*M}\textstyle{\frac 1{n!}} \, \omega_M^n 
\end{equation}
where~$\omega_M$ is the standard symplectic form on~$T^*M$ and~$\frac 1{n!} \, \omega_M^n=\frac 1{n!} \, \omega_M\wedge\dots\wedge\omega_M$ is the corresponding volume form.
Equivalently (see Proposition~\ref{prop:volHTsph}), the Holmes--Thompson volume is given as an integral over the unit sphere bundle by the formula
\begin{equation}\label{eq:volHTsph}
\vol(M) = \frac{1}{\epsilon_nn!} \int_{U^*M} \alpha_M \wedge \omega_M^{n-1}.
\end{equation}
\end{definition}
The factor $\frac 1{\epsilon_n}$ ensures that for Riemannian metrics, the Holmes--Thompson definition of volume agrees with the conventional Riemannian definition.

\section{Integral geometry in Finsler manifolds with almost no trapped geodesics} \label{sec:integral}

The goal of this section is to present versions of two classical formulas in integral geometry, namely the formulas of Blaschke~\cite{blaschke1935integralgeometrie} and Santal\'o~\cite{santalo1952measure,santalo1976integral}, which are in turn generalizations for manifolds of the classical Crofton formulas on the Euclidean plane. 
In~\cite{paiva2006wrong}, Blaschke's formula is proved for Finsler manifolds whose space of geodesics 
is a smooth manifold.
Here, we give slightly more general versions which hold for Finsler manifolds with almost no trapped geodesics (and, in particular, for compact Finsler manifolds with minimizing interior geodesics).
The proofs mimick those given by Blaschke, Santaló, and Álvarez-Paiva--Berck.
However, we give them in full in order to provide additional details and introduce the few extra steps needed for the generalization.

\begin{definition} \label{def:traversing}
Let~$M$ be a Finsler $n$-manifold with nonempty boundary.
A \term{traversing geodesic} of~$M$ is a maximal geodesic $\gamma:[0,\ell(\gamma)] \to M$ which does not intersect~$\partial M$, except at its endpoints where it meets the boundary transversely.
The Finsler manifold~$M$ has \term{almost no trapped geodesics} if for almost every unit tangent vector~$v \in UM$, the maximal geodesic~$\gamma_v$ defined by~$\gamma_v'(0)=v$ reaches the boundary of~$M$ in the future and in the past, that is,~$\gamma_v(t)\in\partial M$ for some~$t\geq 0$ and some~$t\leq 0$.
\end{definition}

For instance, a compact Finsler manifold with minimizing interior geodesics has almost no trapped geodesics. 
Another example is obtained by taking a closed Finsler manifold with ergodic geodesic flow and removing a smoothly bounded nonempty open set.

As we will explain below, the space~$\Gamma$ of traversing geodesics of~$M$ is a~$(2n-2)$-dimensional manifold admitting a natural symplectic structure, whose corresponding natural volume measure is denoted by~$\mu_\Gamma$; see Definition~\ref{def:measure}.

\begin{theorem}[Blaschke's formula] \label{thm:blaschke} 
Let~$M$ be a Finsler $n$-manifold with almost no trapped geodesics. 
Then the Holmes--Thompson volume of an immersed hypersurface~$N \subseteq M$ is equal to
\begin{equation}\label{eq:blaschke}
\vol_{n-1}(N) = \frac{1}{2\,\epsilon_{n-1}} \int_{\gamma\in\Gamma}\#(\gamma \cap N)\,\diff\mu_\Gamma(\gamma)
\end{equation}
where $\#(\gamma \cap N)$ is the number of times that~$\gamma$ intersects~$N$.

Similarly, the Holmes--Thompson volume of a co-oriented immersed hypersurface $N \subseteq M$ is equal to
\begin{equation}\label{eq:blaschke_cooriented}
\vol_{n-1}(N) = \frac{1}{\epsilon_{n-1}}
\int_{\gamma\in\Gamma}\#(\gamma \cap^+ N)\,\diff\mu_\Gamma(\gamma)
\end{equation}
where $\#(\gamma \cap^+ N)$ is the number of times that~$\gamma$ intersects~$N$ transversely in the positive direction.
\end{theorem}

In equation~\eqref{eq:blaschke}, we can restrict the integral to geodesics~$\gamma \in \Gamma$ which are transverse to the hypersurface~$N$ since the geodesics~$\gamma\in\Gamma$ which are tangent to~$N$ form a subset of zero measure; see Proposition~\ref{prop:negligible}.\eqref{negligible3}.

\begin{remark} \label{rem:mu_total}
Since every traversing geodesic intersects~$\partial M$ positively exactly once, we derive from~\eqref{eq:blaschke_cooriented} that the total measure of the space~$\Gamma$ is 
\[  \mu_\Gamma(\Gamma) = \epsilon_{n-1} \vol_{n-1}(\partial M) \]
In particular, if $M$ is compact, then $\mu_\Gamma(\Gamma)<\infty$.
\end{remark}

\begin{theorem}[Santal\'o's formula] \label{thm:santalo} 
Let~$M$ be a Finsler $n$-manifold with almost no trapped geodesics. 
Then the Holmes--Thompson volume of a smoothly-bounded domain~$D\subseteq M$ is equal to
\begin{equation}\label{eq:santalooo}
\vol_n(D) = \frac{1}{n\,\epsilon_n} \int_{\gamma\in\Gamma}\length(\gamma\cap D)\,\diff\mu_\Gamma(\gamma).
\end{equation}
\end{theorem}

In the case of Finsler surfaces with self-reverse metric, 
the Blaschke and Santal\'o formulas specialize as follows.

\begin{corollary}
Let~$M$ be a self-reverse Finsler metric surface with almost no trapped geodesics.
Then the length of any immersed curve~$c$ in~$M$ is
\begin{equation} \label{eq:crofton}
\length(c) = \frac{1}{4} \int_{\gamma\in\Gamma} \#(\gamma \cap c) \, \diff\mu_\Gamma(\gamma),
\end{equation}
and the Holmes--Thompson area of any smoothly-bounded domain~$D \subseteq M$ is
\begin{align}
\area(D)
&= \frac {1}{2\pi} \int_{\gamma\in\Gamma} \length(\gamma \cap D) \, \diff\mu_\Gamma(\gamma) \label{eq:santalo}\\
&= \frac{1}{8\pi} \iint_{(\gamma_0, \gamma_1) \in \Gamma \times \Gamma} \#(\gamma_0 \cap \gamma_1 \cap D) \, \diff\mu_\Gamma(\gamma_0) \, \diff\mu_\Gamma(\gamma_1).\label{eq:croftonsantalo}
\end{align}
\end{corollary}

The equation~\eqref{eq:croftonsantalo}, obtained from~\eqref{eq:santalo} and~\eqref{eq:crofton}, will be called the Santal\'o+Blaschke formula.
In deducing this formula, we use the hypothesis that the metric is self-reverse when we equate the length of a geodesic with its Holmes--Thompson measure.
In general, the Holmes--Thompson measure of a curve is the average of its forward and backward lengths.

The rest of this section is dedicated to describing the symplectic structure on~$\Gamma$ and proving Theorems~\ref{thm:blaschke} and~\ref{thm:santalo}.

\subsection{Symplectic manifold of traversing geodesics}

Let~$M$ be a Finsler $n$-manifold.
Recall that~$\Gamma$ is the space of traversing geodesics of~$M$. 
This space~$\Gamma$ is a $(2n-2)$-dimensional manifold parameterized by the initial vectors $\gamma'(0) \in UM|_{\partial M}$ of the geodesics $\gamma:[0,\ell(\gamma)] \to M$ of~$\Gamma$.
Note that the length~$\ell(\gamma)$ depends smoothly on~$\gamma\in\Gamma$.

\begin{definition} \label{def:measure}
Define the open subset of~$U^*M$
\[ U_\Gamma^*M = \{ \Leg(\gamma'(t)) \in U^*M \mid \gamma\in\Gamma,\ t\in[0,\ell(\gamma)] \} \]
consisting of the momentums of the traversing geodesics of~$M$.
Note that~$U_\Gamma^*M$ is a $Z$-invariant open subset of~$U^*M$.

Consider the surjective submersion 
\begin{equation} \label{eq:submersion}
\pi_\Gamma:U_\Gamma^*M \to \Gamma
\end{equation}
taking any momentum $\xi \in U_\Gamma^*M$ to the geodesic $\gamma \in \Gamma$ that it generates.
The fibers of~$\pi_\Gamma$ are the $Z$-orbits corresponding to the traversing geodesics. 

There exists a unique 2-form ~$\omega_\Gamma$ on~$\Gamma$ such that
\[ \pi_\Gamma^* \, \omega_\Gamma = \omega_M|_{U^*_\Gamma M}. \] 
This follows from the invariance of the 2-form $\omega_M|_{U^*_\Gamma M}$ under the cogeodesic flow, and the fact that this form vanishes in the direction of~$Z$ according to Definition~\ref{def:cogeodesic}. (See also~\cite{cossarini2020finsler} for details, or~\cite[\S4.3]{abraham1978foundations} for a general account on symplectic reduction.)
The form~$\omega_\Gamma$ is symplectic, thus it determines on~$\Gamma$ a smooth volume measure~$\mu_\Gamma$ given by
\begin{equation} \label{eq:muG}
\diff\mu_\Gamma =\frac{1}{(n-1)!} \, |\omega_\Gamma^{n-1}|.
\end{equation}

\end{definition}

\subsection{Non-traversing geodesics are negligible}

We will need the following result in order to establish our versions of Blaschke's and Santal\'o's formulas.
This feature is not required in the previous versions and necessitates the manifold to have almost no trapped geodesics.

\medskip

Recall that a subset~$A$ of a manifold~$X$ is \term{negligible} in~$X$ if the image of~$A$ in any local chart of~$X$ has zero measure.

\begin{proposition} \label{prop:negligible}
\mbox{ }
\begin{enumerate}
\item The complement of the open subset $U^*_\Gamma M \subseteq U^*M$ is negligible in~$U^*M$.\label{negligible1}
\item Given a hypersurface $H \subseteq U^*M$ transverse to~$Z$, the complement of $H \cap U^*_\Gamma M$ is negligible in~$H$. \label{negligible2}
\item The set of geodesics~$\gamma \in \Gamma$ tangent to an immersed hypersurface of~$M$ or passing through an immersed submanifold of~$M$ of codimension~\mbox{$>1$} has zero measure in~$\Gamma$. \label{negligible3}
\end{enumerate}
\end{proposition}

\begin{proof}
To avoid technical problems we extend the Finsler metric to an open manifold~$M^+$; see Remark~\ref{rmk:extend}. This ensures the cogeodesic flow $(\xi,t)\mapsto Z^t(\xi)$ is defined on an open domain.

By definition of~$\Gamma$, the complement~$U^*M \setminus U^*_\Gamma M$ is formed of momentums of two types.
First, momentums of~$U^*M$ that correspond to geodesics of~$M$ with at least one end trapped in~$M$.
These momentums form a negligible set since~$M$ has almost no trapped geodesics; see Definition~\ref{def:traversing}.
Second, momentums of~$U^*M$ corresponding to geodesics tangent to the boundary~$\partial M$. These momentums are of the form~$Z^t(\xi)$, where~$Z^t$ is the cogeodesic flow and~$\xi$ is the Legendre image of a unit vector~$v$ tangent to~$\partial M$. These unit vectors form a manifold~$U\partial M$ of dimension~$2n-3$.
Thus, by Sard's theorem, the map from an open subset of $U\partial M \times \R$ to~$U^*M$ defined by $(v,t) \mapsto Z^t(\Leg(v))$ has negligible image in~$U^*M$.
Having considered both types of momentums, we conclude that the complement~$U^*M \setminus U^*_\Gamma M$ is negligible in~$U^*M$.

For the second point, simply observe that if~$A$ is a~$Z$-invariant negligible subset of~$U^*M$ and~$H$ is a hypersurface of~$U^*M$ transverse to~$Z$, then~$A \cap H$ is negligible in~$H$.
Apply this property to~$A=U^*M \setminus U^*_\Gamma M$ to conclude.

The proof of the third point is similar to the proof of the first point and relies on Sard's theorem.
\end{proof}

\subsection{Manifold of positive momentums across a hypersurface}

We will need the following notion in the proof of Blaschke's formula.

\begin{definition}
Let~$N$ be a co-oriented embedded hypersurface in a Finsler manifold~$M$.
(For example, we can have~$N=\partial M$ co-oriented so that inwards-pointing vectors are positive.)
Denote by~$C^*N \subseteq U^*M|_{N}$ the \term{manifold of momentums crossing positively} the hypersurface~$N$, that is, the momentums corresponding under the Legendre map to unit vectors transverse to~$N$ pointing in the positive direction according to the co-orientation of~$N$.
Note that~$C^*N$ is an open subset of~$U^*M|_{N}$ and therefore a differentiable manifold, with the structure of an open ball bundle over~$N$.

Consider the \term{restriction map}
\[\begin{array}{rccc}
\rho_N:
&C^*N &\longrightarrow&\Int(B^*N) \\
& \xi \in T_x^*M &\longmapsto    & \xi|_{T_xN}
\end{array}\]
to the interior~$\Int(B^*N)$ of the unit co-ball bundle~$B^*N$ of~$N$.
\end{definition}

The following statement can be found in~\cite[Lemma~5.4]{paiva2006wrong}.
We simply provide the details of the proof.

\begin{lemma} \label{lem:symplectomorphism} 
The space~$C^*N$ is a symplectic submanifold of~$T^*M$ and the restriction map
\[ \rho_N: (C^*N,\omega_M) \to (\Int(B^*N),\omega_N) \]
is a symplectomorphism.
Thus,
\[ \rho_N^* \, \omega_N = {\omega_M}|_{C^*N}. \]
\end{lemma}

\begin{proof} Let $\xi \in C^*N$ with basepoint~$x \in N$.
By definition, the norm of~$\xi$ is~$1$, so the norm of its restriction~$\xi'$ to~$T_xN$ is at most~$1$.
Furthermore, by strong convexity of~$F_x^*$, the linear form~$\xi$ attains its maximum only at its Legendre-dual unit vector, which is positive and thus not contained in~$T_xN$.
Therefore, $\| \xi' \| <1$ and the restriction map~$\rho_N$ takes values in~$\Int(B^*N)$.

To see that~$\rho_N$ is a diffeomorphism, we employ local coordinates $(x_i)_{1\leq i\leq n}$ in~$M$ so that the hypersurface~$N$ is given by the equation~$x_n=0$. Let $(x_i,v_i)_i$ and $(x_i,\xi_i)_i$ be the corresponding coordinates in~$TM$ and~$T^*M$. In terms of these coordinates, the operator~$\rho_N$ acts by supressing the last coefficient, that is, if~$\xi=(\xi_i)_{1\leq i\leq n}$, then~$\xi'=(\xi_i)_{1\leq i\leq n-1}$. Hence~$\rho_N$ is smooth.

To prove that~$\rho_N$ is bijective, consider a covector~$\xi'=(\xi_i)_{1\leq i\leq n-1}\in \Int(B_x^*N)$ and denote its norm~$\lambda=\|\xi'\|<1$. The covectors~$\xi\in T^*_xM$ such that~$\xi|_{T_xN}=\xi'$ are of the form~$\xi^t=(\xi_1,\dots,\xi_{n-1},t)$ with~$t\in\R$.
Consider the function $t\mapsto \|\xi^t\|$, where~$\|\cdot\|$ is the norm~$F_x^*$ on~$T_x^*M$ that is dual to~$F_x$.
This function is bounded below by~$\lambda$, and by the Hahn--Banach theorem, this lower bound is attained at some~$t_0\in\R$.
Furthermore, since the norm~$F_x^*$ is strongly convex, the set of values of~$t$ such that~$\|\xi^t\|\leq 1$ is a compact interval~$[t_-,t_+]$ that contains~$t_0$ in its interior, and~$\|\xi^t\|=1$ if and only if~$t=t_\pm$.
Thus we are left with two candidates~$\xi^{t_\pm}$ that are the only unit covectors~$\xi$ whose restriction to~$T_xN$ is~$\xi'$.

We claim that~$\xi^{t_+}$ is positive (and~$\xi^{t_-}$ is negative).
That is, the vector that is in Legendre correspondence with $\xi^{t_+}$ (\ie, the unit vector where $\xi^{t_+}$ attains its norm) 
is positive.
Indeed, when~$t=t_0$, the covector~$\xi^t$, as a function~$B_xM\to\R$, is bounded above by~$\lambda$.
As~$t$ increases towards~$t_+$, the coefficient~$\xi_n$ increases, and thus the values of~$\xi^t(v)$ for $v$ on the 
negative side decrease (hence they are~$<\lambda$).
Thus, any functional~$\xi^t$ with~$t>t_0$, restricted to the ball~$B_xM$, must attain its maximum value~$\|\xi^t\|$ (which is $>\lambda$) on a positive vector, as required.
This shows that~$\xi^t$ is positive if $t>t_0$ (and, similarly, $\xi^t$ is negative if~$t<t_0$).
We conclude that~$\xi^{t_+}$ is the only positive unit covector~$\xi$ whose restriction to~$T_xN$ is~$\xi'$. 
This proves that~$\rho_N$ is bijective.
Additionally,~$t_+$ depends smoothly on~$\xi'$ by the implicit function theorem. 
This finishes the proof that the restriction map $\rho_N:C^*N\to\Int(B^*N)$ is a diffeomorphism.

Let us show that~${\rho_N}^*\alpha_N = \alpha_M|_{C^*N}$.
In canonical coordinates, the tautological 1-form~$\alpha_M$ on~$T^*M$ is written as~$\alpha_M = \sum_{i=1}^n \xi_i \diff x_i$.
In restricting to~$C^*N$, the last term vanishes because~$x_n=0$ on~$N$, thus the restricted form can be written as ${\alpha_M}|_{C^*N} = \sum_{i=1}^{n-1} \xi_i \diff x_i$.
On the other hand, the tautological 1-form of~$N$ is $\alpha_N = \sum_{i=1}^{n-1} \xi_i \, \diff x_i$, and this expression is unchanged by the pullback~${\rho_N}^*$ since the map~$\rho_N:C^*N\to \Int(B^*N)$ acts simply by suppressing the coordinate~$\xi_n$.
We conclude that~${\rho_N}^*\alpha_N = {\alpha_M}|_{C^*N}$.
Taking the exterior differential of this expression, we obtain~${\rho_N}^*\omega_N = \omega_M|_{C^*N}$.
This implies that~$C^*N$ is a symplectic submanifold of~$T^*M$.
%
\end{proof}

\subsection{Coarea formula and fiber integration}

In the proofs of Blaschke's and Santal\'o's formulas,
we will need the following version of the coarea formula; see~\cite[(16.24.8)]{dieudonne1972treatise} (see also \cite[Theorem~3.2.3]{federer1969geometric} and~\cite[Theorem~5.5.8]{burago2001course} when~$n=m$). 

\begin{lemma} \label{lem:coarea}
Let $\pi:X \to Y$ be a submersion between two oriented manifolds of dimension~$n$ and~$m$ with~$n \geq m$.
Let~$\alpha$ and~$\beta$ be two differential forms on~$X$ and~$Y$ of degree~$n-m$ and~$m$. Then
\[  \int_X \alpha \wedge \pi^*\beta
  = \int_{y \in Y} \left( \int_{\pi^{-1}(y)}\alpha \right)\,\beta \]
where~$\pi^{-1}(y)$ is endowed with the orientation induced by~$\pi$ from the orientations of~$X$ and~$Y$.

In particular, for~$n=m$ and $\alpha = 1$, we have
\begin{equation} \label{eq:coarea n=m}
\int_X \pi^*\beta = \int_{y \in Y} \# (\pi^{-1}(y))\,\beta.
\end{equation}
\end{lemma}

\subsection{Proof of the Blaschke formula}

We can now proceed to the proof of Blaschke's formula~\eqref{eq:blaschke}.

\begin{proof}[Proof of Theorem~\ref{thm:blaschke}]
We will follow the proof given in~\cite[Theorem~5.2]{paiva2006wrong} under the extra assumption that the space of oriented geodesics on~$M$ is a manifold.

The Blaschke formula~\eqref{eq:blaschke} for a non-cooriented hypersurface~$N$ can be deduced from the co-oriented version~\eqref{eq:blaschke_cooriented} by taking the co-oriented double cover of~$N$. 
Therefore it is sufficient to prove the latter formula.
Furthermore, every immersed hypersurface can be decomposed into a disjoint union of embedded hypersurfaces up to a negligible set. 
Therefore it is sufficient to prove~\eqref{eq:blaschke_cooriented} for a co-oriented embedded hypersurface~$N$.

By definition of the Holmes--Thompson volume, see~\eqref{eq:volHT}, we have
\begin{align*}
\vol_{n-1}(N)
& =  \frac{1}{\epsilon_{n-1}(n-1)!} \int_{B^*N} \omega_N^{n-1} \\
& = \frac{1}{\epsilon_{n-1}(n-1)!} \int_{C^*N} \omega_M^{n-1}
\end{align*}
where the second equality follows from Lemma~\ref{lem:symplectomorphism}.

Now, apply Proposition~\ref{prop:negligible}.\eqref{negligible2} with $H=C^*N \subseteq U^*M$.
It follows that $C^*N \cap U_\Gamma^*M$ has full measure in~$C^*N$.
Thus,
\[  \int_{C^*N} \omega_M^{n-1}
= \int_{C^*N \, \cap \, U_\Gamma^*M} \omega_M^{n-1}.  \]

Consider the map \mbox{$\pi:C^*N \cap U_\Gamma^*M \to \Gamma$} taking a unit momentum of~$M$ based at~$N$ pointing in a positive direction (with respect to the co-orientation of~$N$) to the traversing geodesic it generates.
Apply the fiber integration formula~\eqref{eq:coarea n=m} to this map with~$\beta= \omega_\Gamma^{n-1}$.
This yields the relation
\[ \int_{C^*N \, \cap \, U_\Gamma^*M} \omega_M^{n-1}
= \int_{\gamma \in \Gamma} \#(\gamma \cap^+ N) \, \omega_\Gamma^{n-1}. \]
where $\#(\gamma \cap^+N)$ is the number of times that~$\gamma$ crosses~$N$ transversely in the positive sense (as determined by the co-orientation of~$N$).
Taking into account the definition of~$\mu_\Gamma$ by equation~\eqref{eq:muG}, Blaschke's formula follows.
\end{proof}

\subsection{Proof of the Santal\'o formula}

We will need the following proposition expressing the Holmes--Thompson volume of a manifold as an integral over the bundle of dual unit spheres (instead of dual unit balls). 

\begin{proposition} \label{prop:volHTsph}
The Holmes--Thompson volume of a Finsler $n$-manifold~$M$ is equal to
\[ \vol(M)
= \frac{1}{\epsilon_nn!} \int_{U^*M} \alpha_M \wedge \omega_M^{n-1}.\]
\end{proposition}

\begin{proof}
We may assume that $M$ is a compact manifold with corners.
(If $M$ is not compact, we can triangulate it and apply the proposition on each $n$-simplex to infer that it holds on the whole manifold.)

Observe that $\diff(\alpha_M \wedge \omega_M^{n-1}) = \omega_M^n$ and $\partial B^*M = U^*M \cup B^*M|_{\partial M}$.
By Stokes' theorem, we have
\begin{align} \vol_n(M)
& = \frac{1}{\epsilon_n n!} \int_{B^*M} \diff(\alpha_M \wedge \omega_M^{n-1})
\nonumber \\
& = \frac{1}{\epsilon_n n!} \int_{U^*M} \alpha_M \wedge \omega_M^{n-1}
+ \frac{1}{\epsilon_n n!} \int_{B^*M|_{\partial M}} \alpha_M \wedge \omega_M^{n-1},
\nonumber
\end{align}
where $\partial M$ is considered as a piecewise-smooth $(n-1)$-manifold (and we may restrict the integral to its smooth part).


To finish the proof, we shall show that the $(2n-1)$-form $\alpha_M \wedge \omega_M^{n-1}$ vanishes on~$B^*M|_{\partial M}$, hence the second term 
vanishes.
Recall that the tangent space to $B^*M|_{\partial M}$ at~$(x,\xi)$ decomposes as
\[ T_{(x,\xi)} B^*M|_{\partial M} \simeq T_x \partial M \oplus T_x^*M \]
where the horizontal space~$T_x \partial M$ is of dimension~$n-1$ and the vertical space~$T_x^*M$ is of dimension~$n$.
Note that the one-form~$\alpha_M$ vanishes on the vertical space~$T_x^*M$ and the two-form~$\omega_M$ vanishes at bi-vectors formed of two horizontal or two vertical vectors. This follows from the coordinate expression~\eqref{eq:alpha_omega_canon}.
The $(2n-1)$-form~$\alpha_M \wedge \omega_M^{n-1}$ evaluated at $(u_1,\dots,u_{2n-1})$, where~$n-1$ vectors~$u_i$ are horizontal and~$n$ vectors~$u_i$ are vertical, can be written as a sum of terms of the form 
\begin{equation} \label{eq:term}
\pm \alpha_M(u_{\sigma(1)}) \cdot \omega_M(u_{\sigma(2)}, u_{\sigma(3)}) \cdot \ldots \cdot \omega_M(u_{\sigma(2n-2)}, u_{\sigma(2n-1)})
\end{equation}
where~$\sigma$ is a permutation.
If~$u_{\sigma(1)}$ is vertical then the factor~$\alpha_M(u_{\sigma(1)})$ is equal to zero.
If~$u_{\sigma(1)}$ is horizontal then there are only~$n-2$ horizontal vectors (and~$n$ vertical ones) among the remaining vectors, which implies that one of the factors $\omega_M(u_{\sigma(2k)}, u_{\sigma(2k+1)})$ has two horizontal vectors and therefore vanishes.
In both cases, the term~\eqref{eq:term} vanishes.
\end{proof}

Let us prove Santal\'o's formula~\eqref{eq:santalooo}.

\begin{proof}[Proof of Theorem~\ref{thm:santalo}]
Recall that $\omega_M = \pi^*_\Gamma \, \omega_\Gamma$, see Definition~\ref{def:measure}, and that~$U^*_\Gamma M$ has full measure in~$U^* M$, see Proposition~\ref{prop:negligible}.\eqref{negligible1}.
By 
Proposition~\ref{prop:volHTsph}
we have
\[ \vol_n(D)
= \frac{1}{\epsilon_n n!} \int_{U^*D \, \cap \, U^*_\Gamma M} \alpha_M \wedge \pi^*_\Gamma \, \omega_\Gamma^{n-1}. \]
By Lemma~\ref{lem:coarea}, integrating along the fibers of the submersion $\pi:U^*D \cap U^*_\Gamma M \to \Gamma$ induced by~$\pi_\Gamma$, see~\eqref{eq:submersion}, we obtain 
\[ \vol_n(D)
= \frac{1}{\epsilon_n\,n!} \int_{\gamma \in \Gamma} \left( \int_{\pi^{-1}(\gamma)} \alpha_M \right) \, \omega_\Gamma^{n-1} . \]
Since all the fibers~$\pi^{-1}(\gamma) = \{ \Leg(\gamma'(t)) \in U^*M \mid t \in [0,\ell(\gamma)] \} \cap U^*D$ are tangent to the cogeodesic vector field~$Z$ on~$U^*M$ and~$\alpha_M(Z) =1$, we derive
\[ \int_{\pi^{-1}(\gamma)} \alpha_M = \length(\gamma \cap D). \]
Hence,
\begin{align*}
\vol_n(D)
& = \frac{1}{\epsilon_n\,n!}
\int_{\gamma \in \Gamma} \length(\gamma \cap D) \, \omega_\Gamma^{n-1} \\
& = \frac{1}{\epsilon_n\,n}
\int_{\gamma\in\Gamma}\length(\gamma\cap D)\,\diff\mu_\Gamma(\gamma).
\end{align*}
\end{proof}

\section{Discretization of Finsler surfaces} \label{sec:discretization}

The goal of this section is to describe a discretization of Finsler disks with minimizing interior geodesics into simple discrete metric disks.
For this, we adapt the general approach of discretization developed in~\cite{cossarini2018discrete} in relation with the filling area conjecture.
The main novelty is that, in our case, the discrete geometry is described by a system of curves (wall system) made of geodesics.

\medskip

First, we need to fix some notation regarding intersections of maps.

\begin{definition}\label{def:map_intersections}
The \term{intersections} of a map $f:X\to Y$ with a map $f':X'\to Y$ lying in a subset $A\subseteq Y$ are the ordered pairs in the set
\[ I_A(f,f') = \left\{ (x,x')\in X\times X'
\text{ such that } f(x) = f'(x') \in A \right\}. \]
The \term{number of intersections} between~$f$ and~$f'$ is defined as
\[ \# (f\cap f') = \# I_Y(f,f'),\]
where~$\#S$ denotes the cardinality of a set~$S$.

Similarly, the \term{self-intersections} of a map $f:X\to Y$ lying in a subset~$A\subseteq Y$ are the \term{unordered} pairs in the set
\[ I_A(f) = \left\{ \{x,x'\} \subseteq X \text{ such that }
f(x) = f(x') \in A \text{ but } x\neq x' \right\}, \]
and the \term{multiplicity} of a point~$y\in Y$ as a self-intersection of~$f$ is the number~$\# I_{\{y\}}(f)$.
A self-intersection is \term{simple} if it has multiplicity 1.
\end{definition}

Let us introduce the notion of wall system on a disk; see~\cite{cossarini2018discrete}.

\begin{definition} \label{def:wall}
A (smooth) \term{wall system} on a surface~$M$ is a $1$-dimensional (smooth) immersed submanifold~$\W$ satisfying the following conditions:
\begin{enumerate}
\item the immersion map is proper (that is, the preimage of any compact subset of~$M$ is compact);
\item $\W$ is transverse to the boundary~$\partial M$ and satisfies $\partial\W = \W \cap \partial M$;
\item $\W$ is self-transverse and has only simple self-intersections;
\item\label{item:bound} no self-intersections of~$\W$ lie on the boundary~$\partial M$.
\end{enumerate}

As a technical remark, we note that the symbol~$\W$ denotes the \term{immersion map}, not its image $\Img(\W)\subseteq M$, nor its domain. 
The domain is a $1$-manifold, \ie, a disjoint union of countably many intervals and circles. 
Hence the expression $\partial\W\subseteq\partial M$ involves an abuse of notation and actually means $\Img(\partial\W)\subseteq\partial M$, where~$\partial \W$ is the restriction of the map~$\W$ to the boundary of the domain of~$\W$.
The image of~$\W$ will also be denoted~$\W$.
Thus, the expression~$M\setminus\W$ denotes $M\setminus\Img(\W)$.

Eventually we will need to relax the definition by dropping condition~\eqref{item:bound}. In this case, we say that~$\W$ is a \term{quasi wall system} on~$M$.

The curves that form a (quasi) wall system are called its \term{walls}. Note that if the surface~$M$ is compact, then~$W$ consists of finitely many compact walls; each of these walls is either a loop that avoids the boundary or an arc that meets the boundary only at its two endpoints.

A quasi wall system~$\W$ on a disk~$D$ is \term{simple} if its walls are arcs that have no self-intersections and that meet each other at most once.\footnote{Simple wall systems are also called pseudoline arrangements; see~\cite{cossarini2018discrete}. However, some authors (\eg,~\cite{felsner2018pseudoline}) only consider complete pseudoline arrangements, which are those where every pair of walls crosses \emph{exactly} once.}

In this paper, every quasi wall system~$\W$ is smooth unless we make it clear that it is piecewise smooth. In that case, the non-smooth points of~$\W$ may not coincide with the self-intersection points of~$\W$.
 Note that a piecewise smooth quasi wall system can be turned into a smooth quasi wall system by an isotopic deformation.
\end{definition}

\begin{example}
Let~$D$ be the unit disk in the Euclidean plane.
A wall system made of the horizontal and vertical diameters of~$D$ has area~$1$.
A quasi wall system made of the three sides of an inscribed triangle of~$D$ has area~$\frac{3}{2}$.
\end{example}

We will also need the following definitions regarding the geometry induced by a quasi wall system.

\begin{definition}
Every quasi wall system~$\W$ on a compact surface~$M$ determines a \term{discrete length}
\begin{equation} \label{eq:lengthW}
\length_\W(c) = \#(c \cap \W)
\end{equation}
for curves~$c$ in~$M$.
That is, the length of a curve is the number of times it intersects the quasi wall system (counted with multiplicity).
Every quasi wall system~$\W$ also induces a pseudo-distance on~$M \setminus\W$ defined by
\begin{equation*} 
d_\W(x,y) = \inf_c \length_\W(c)
\end{equation*}
where the infimum is taken over all paths of~$M$ joining~$x$ to~$y$.
We will refer to the pseudo-distance~$d_\W$ on~$M$ as the \term{discrete distance} induced by~$\W$ on~$D$.

The \term{discrete area} of~$(M,\W)$ is the number of self-crossings of~$\W$ contained in the interior of~$M$ plus half the number of self-crossings on the boundary.
That is,
\[ \area(M,\W)
= \# I_{\Int M}(\W) + \tfrac 12 \# I_{\partial M}(\W)
= \# I_{\Int M + \frac 12 \partial M}(\W)   \]
where $\# I_{\Int M + \frac 12 \partial M}$ is just an abbreviation for $\# I_{\Int M} + \tfrac 12 \# I_{\partial M}$.
\end{definition}
Note that, if~$\W$ consists of finitely many curves~$\gamma_i$, then
\begin{equation}
\area(M,\W)
= \sum_{i<j}\# I_{\Int M + \frac 12 \partial M}(\gamma_i,\gamma_j)
+ \sum_{i}\# I_{\Int M + \frac 12 \partial M}(\gamma_i) \label{eq:area}
\end{equation}
When the quasi wall system is simple, the curves of~$\W$ have no self-intersections and the second sum vanishes.

\medskip

We will need the following result describing the intersection of two distance-realizing arcs of~$M$.
Recall that~$\Gamma$ is the space of traversing geodesics of~$M$ (\ie, geodesic arcs of~$M$ which do not intersect~$\partial M$ except at their endpoints, where they meet the boundary transversely).

\begin{lemma} \label{lem:intersection}
Let~$M$ be a self-reverse Finsler metric disk with minimizing interior geodesics.
Let~$\gamma\in\Gamma$ be a traversing geodesic of~$M$ and let~$[x,y]$ be a distance-realizing arc of~$M$ with endpoints~$x$ and~$y$ not lying in~$\gamma$.
Then
\[ \#(\gamma \cap [x,y]) =
\begin{cases}
1 & \text{if } \gamma \text{ separates } x \text{ and } y \\
0 & \text{otherwise}
\end{cases}
\]
\end{lemma}

\begin{proof}
By Theorem~\ref{thm:shortest_are_C1}, the distance-realizing arc~$[x,y]$ is~$C^1$. 

Suppose that the arcs~$\gamma$ and~$[x,y]$ are tangent, either at an interior point of~$M$ or at an endpoint of~$\gamma$ in~$\partial M$.
In both cases, this implies that~$[x,y]$ contains~$\gamma$ since the distance-realizing arc~$[x,y]$ follows the geodesic flow in the interior of~$M$ and the endpoints~$x,\,y$ do not lie in~$\gamma$.
Now, since the interior geodesic~$\gamma$ is transverse to~$\partial M$ at its endpoints~$\bar{x}$ and~$\bar{y}$, the distance-realizing arc~$[x,y]$ is not differentiable at~$\bar{x}$ and~$\bar{y}$.
In particular, it is not~$C^1$, which is absurd.
Therefore, the arcs~$\gamma$ and~$[x,y]$ may only have transverse intersections.

Suppose that the arcs~$\gamma$ and~$[x,y]$ intersect at least twice, say at~$a$ and~$b$ (with~$a$ and~$b$ different from~$x$ and~$y$).
Since both arcs are distance-realizing curves, the subarcs~$[a,b] \subseteq [x,y]$ and $\gamma_{ab} \subseteq \gamma$ joining~$a$ and~$b$ have the same length.
Construct an arc~$\alpha$ joining~$x$ and~$y$ by replacing the subarc~$[a,b]$ of~$[x,y]$ with the arc~$\gamma_{ab}$ of the same length.
By construction, the arc~$\alpha$ is a distance-realizing curve.
But since the intersection between~$\gamma$ and~$[x,y]$ is transverse, the arc~$\alpha$ is not differentiable at~$a$ and~$b$.
In particular, it is not~$C^1$, which is absurd.
Therefore, the arcs~$\gamma$ and~$[x,y]$ intersect at most once, and so exactly once if~$\gamma$ separates~$x$ and~$y$.

Suppose now that~$\gamma$ does not separate~$x$ and~$y$.
Then the arc~$[x,y]$ does not intersect~$\gamma$.
Otherwise, it would go from one side of~$\gamma$ to the other (recall that~$\gamma$ and~$[x,y]$ have transverse intersection) and, because~$x$ and~$y$ are on the same side of~$\gamma$, it would have to cross~$\gamma$ a second time, which is excluded.
Therefore, the arcs~$\gamma$ and~$[x,y]$ do no intersect if~$\gamma$ does not separate~$x$ and~$y$.
\end{proof}

Let us compare the shortest paths for Finsler metrics and discrete metrics.

\begin{definition}
A quasi wall system is \term{geodesic} if its walls are geodesics.
\end{definition}

\begin{proposition} \label{prop:min}
Let~$M$ be a self-reverse Finsler metric disk with minimizing interior geodesics, and let~$\W$ be a geodesic quasi wall system on~$M$.
Then, every distance-realizing arc~$[x,y]$ of~$M$ with endpoints~$x,\,y$ not lying in~$\W$ is also length minimizing with respect to~$\W$.\\
Thus, for every $x,y \in M \setminus\W$, we have
\begin{equation} \label{eq:min}
d_\W(x,y) = \length_\W([x,y]).
\end{equation}
\end{proposition}

\begin{proof}
The quasi wall system~$\W$ is made of finitely many geodesics~$\gamma_i$ that are transverse to~$\partial M$.
By Lemma~\ref{lem:intersection}, the arc~$[x,y]$ crosses only those geodesics~$\gamma_i$ that separate~$x$ from~$y$, exactly once. Therefore, no curve from~$x$ to~$y$ can be shorter than~$[x,y]$ with respect to~$\W$.
\end{proof}

Before proceeding we derive a useful consequence of the last lemma.

\begin{lemma}\label{thm:semiper} Let~$M$ be a self-reverse Finsler metric disk with minimizing interior geodesics. Then
\[ d(x,y) \leq \tfrac 12 \length(\partial M) \]
for any pair of points~$x,y\in M$. The same inequality holds if the distance and length are taken with respect to a geodesic quasi wall system~$\W$, that is,
\[ d_\W(x,y) \leq \tfrac 12 \length_\W(\partial M) \]
for points $x,y\in M\setminus\W$.
\end{lemma}

\begin{proof} Join the points $x,y\in M$ by a distance-realizing arc $[x,y]$. By Lemma~\ref{lem:intersection}, each traversing geodesic~$\gamma$ of~$M$ intersects~$[x,y]$ at most once and meets~$\partial M$ exactly twice. Then the inequality $d(x,y)\leq\frac 12\length(\partial M)$ follows from Blaschke's formula~\eqref{eq:crofton} applied to~$[x,y]$.

The claim regarding the geodesic quasi wall system~$\W$ is proved in a similar way.
By Proposition~\ref{prop:min}, the distance-realizing arc~$[x,y]$ is also length-minimizing with respect to~$\W$.
Since each wall of~$\W$ crosses~$[x,y]$ at most once and meets~$\partial M$ exactly twice, we derive the desired second inequality from the definition of~$\length_\W$; see~\eqref{eq:lengthW}.
\end{proof}

Simple wall systems can be used to discretize Finsler disks~$M$ with minimizing interior geodesics.

\medskip

For every $a,b \in \R$ and every~$\eps >0$, we write $a \simeq b \pm \eps$ if $|a-b| < \eps$.

\begin{theorem} \label{theo:existence}
Let~$(M,F)$ be a self-reverse Finsler metric disk with minimizing interior geodesics.
Then, for every~$\eps >0$ and every integer~$n$ large enough, there exists a wall system~$\W$, made of~$n$ geodesics of~$M$, such that for every $x,y \in M \setminus\W$, we have
\begin{equation} \label{eq:distW}
\frac{1}{n} d_\W(x,y) \simeq  \frac{2}{L} d_F(x,y) \pm \eps
\end{equation}
\begin{equation} \label{eq:areaW}
\frac{2}{n^2-n} \, \area(M,\W) \simeq \frac{2\pi}{L^2} \, \area(M,F) \pm \eps
\end{equation}
where $L=\length_F(\partial M)$.
Furthermore, the wall system~$\W$ is necessarily simple.
\end{theorem}

Note that \cite[Theorem 7.1]{cossarini2018discrete} states the existence of a theorem with similar approximation properties but not necessarily made of geodesics.

\begin{proof}
The wall system~$\W$ will be made of random geodesics.
Recall that~$\Gamma$ is the space of traversing geodesics of~$M$ (\ie, geodesic arcs of~$M$ which do not intersect~$\partial M$ except at their endpoints where they meet the boundary transversely) and has a natural measure~$\mu_\Gamma$; see~\eqref{eq:muG}. Furthermore, this space has finite total measure~$\mu_\Gamma(\Gamma) = 2L$; see~\ref{rem:mu_total}. Thus we may define on~$\Gamma$ the probability measure~$\Prob=\frac{\mu_\Gamma}{2L}$.

Take~$n$ independent identically distributed (i.i.d.) random geodesics~$\gamma_1$, \dots, $\gamma_n$ of~$\Gamma$
with probability distribution~$\Prob$.
Almost surely, these geodesics form a wall system~$\W$ of~$M$; see Definition~\ref{def:wall}; because they are pairwise different and form only simple crossings located in the interior of~$M$. Moreover, this wall system is simple, since the geodesics are minimizing and therefore they cannot cross each other more than once by Lemma~\ref{lem:intersection}.
At this point, Theorem~\ref{theo:existence} follows from the next two lemmas.

\medskip

The first lemma is obtained by applying the weak law of large numbers to the Blaschke formula~\eqref{eq:crofton} in a uniform way.

\begin{lemma} \label{lem:dd}
With probability converging to 1 as~$n\to \infty$, the estimate
\[ \frac{1}{n} d_\W(x,y) \simeq  \frac{2}{L} d_F(x,y) \pm \eps \]
holds for every $x,y \in M \setminus \W$.
\end{lemma}

\begin{proof}
Let~$\DD$ be a finite covering of~$M$ by smoothly bounded disks~$D$ with perimeter $\length_F(\partial D) < \eps$.
Fix a basepoint~$p$ in each disk~$D \in \DD$ and denote by~$P$ the collection of all basepoints.
Almost surely, the geodesics of~$\W$ avoid the points of~$P$ and are transverse to the boundaries of the disks~$D\in\DD$.

The following claim shows that the conclusion of the lemma holds in some finite cases.

\begin{claim}
The following assertions hold with probability converging to 1 as~$n\to\infty$.
\begin{enumerate}
\item For every pair of points $p,q\in P$, we have
\begin{equation} \label{eq:claim1}
\frac{1}{n} d_\W(p,q) \simeq  \frac{2}{L} d_F(p,q) \pm \eps.
\end{equation} \label{item1}
\item For every disk~$D\in\DD$ and every pair of points $x,y \in D \setminus \W$, we have
\begin{equation} \label{eq:diam}
\frac{1}{n} d_\W(x,y) \leq \left( \frac{1}{L} + \frac{1}{2} \right) \varepsilon.
\end{equation}
\label{item2}
\end{enumerate}
\end{claim}

\begin{proof}
\eqref{item1} 
Recall that the distance-realizing arc~$[p,q]$ is~$C^1$ embedded in~$M$; see Theorem~\ref{thm:shortest_are_C1}.

The intersection function $f=f_{p,q}:\Gamma \to \N$ defined by
\[ f(\gamma) = \#(\gamma \cap [p,q]) \]
is a nonnegative measurable function.
By Blaschke's formula~\eqref{eq:crofton}, the random variables $X_i=f(\gamma_i)$ with $1 \leq i \leq n$ are i.i.d. with finite expected value
\[ \Expe(X_i) 
= \int_\Gamma \#(\gamma \cap [p,q]) \, \dProb = \frac{2}{L} d_F(p,q). \]
Note that $\Expe|X_i| = \Expe(X_i) < \infty$.
By the weak law of large numbers applied to~$\{X_i\}$ (see e.g. ~\cite{tao2008strong}), we derive that
\[ \left| \frac{1}{n}
\sum_{i=1}^n \#(\gamma_i \cap [p,q]) - \frac{2}{L} d_F(p,q) \right|
< \eps \]
with probability converging to 1 as~$n \to \infty$.
By Proposition~\ref{prop:min}, we have
\[d_\W(p,q) = \length_\W([p,q])=\sum_{i=1}^n \#(\gamma_i \cap [p,q]),\]
hence \eqref{item1} follows.
\medskip

\eqref{item2} The proof of the second assertion is similar.
For a disk~$D\in\DD$, the intersection function $f(\gamma) = \#(\gamma \cap \partial D)$ has expected value $\frac 2L\length_F(\partial D)$ by Blaschke's formula~\eqref{eq:crofton}.
Applying the weak law of large numbers to the random variables $X_i=f(\gamma_i)$ as previously, we derive
\[ \left|
\frac{1}{n} \sum_{i=1}^n \#(\gamma_i \cap \partial D)
- \frac{2}{L} \length_F(\partial D)
\right| < \eps \]
with probability converging to 1 as~$n \to \infty$.
Thus, 
\begin{align*}
\frac{1}{n} \length_{\W}(D)
 &\simeq \frac{2}{L} \length_F(D) \pm \eps \nonumber\\
 &\leq \left(\frac 2L + 1\right) \eps.
\end{align*} 
Since~$D$ is a disk with minimizing interior geodesics, the discrete part of Lemma~\ref{thm:semiper} yields~\eqref{item2}.
\end{proof}

Without loss of generality, we can assume that the conclusion of the previous claim is satisfied.
Let $x,y \in M \setminus\W$.
The points~$x$ and~$y$ lie in some disks~$D_x$ and~$D_y$ of~$\DD$.
Denote by~$p_x$ and~$p_y$ the basepoints of~$D_x$ and~$D_y$.
Since~$D_x$ is a disk with minimizing interior geodesics, by Lemma~\ref{thm:semiper} we have
\[ d_F(x,p_x)\leq \frac12 \length_F(\partial D_x) < \frac 12\eps, \]
thus by the triangle inequality, we obtain
\begin{equation} \label{eq:dxypxpy}
\left| d_F(x,y) - d_F(p_x,p_y) \right| \leq d_F(x,p_x) + d_F(y,p_y) < \eps.
\end{equation}
Combining the triangle inequality with~\eqref{eq:diam}, we obtain
\begin{align} 
\left| \tfrac 1n d_\W(x,y) - \tfrac 1n d_\W(p_x,p_y) \right|
&\leq \tfrac 1n d_\W(x,p_x) + \tfrac 1n d_\W(y,p_y)\nonumber\\
&\leq \left(\frac 2L + 1\right) \eps.\label{eq:dWnxypxpy}
\end{align}
Thus, the following equalities
\[ \tfrac{1}{n} \, d_\W(x,y) \underset{\eqref{eq:dWnxypxpy}}
{\simeq} \tfrac{1}{n} \, d_\W(p_x,p_y) \underset{\eqref{eq:claim1}}
{\simeq} \frac{2}{L} \, d_F(p_x,p_y) \underset{\eqref{eq:dxypxpy}}
{\simeq} \frac{2}{L} \, d_F(x,y) \]
hold up to additive constants which are universal multiples of~$\eps$ (namely, $\left(\frac 2L + 1\right)\eps$ for the first one,~$\eps$ for the second and~$\frac{2}{L}\eps$ for the third one).
Therefore,
\[
\left| \tfrac{1}{n} \, d_\W(x,y) - \frac{2}{L} \, d_F(x,y) \right| < C_0 \, \eps
\]
where $C_0 = \frac{4}{L}+2$.
Hence the first lemma.
\end{proof}

The second lemma is obtained by applying a (slightly generalized) weak law of large numbers to the Santal\'o+Blaschke formula~\eqref{eq:croftonsantalo}.

\begin{lemma} \label{lem:area}
With probability converging to 1 as~$n \to \infty$, we have
\[
\frac{2}{n^2-n} \, \area(M,\W)
\simeq \frac{2\pi}{L^2} \,\area(M) \pm \eps.
\]
\end{lemma}


\begin{proof}
The intersection counting function $f:\Gamma \times \Gamma \to \N$ defined by
\[ f(\gamma,\gamma') = \#(\gamma \cap \gamma') \]
is a measurable function that takes value 0 or 1 almost surely.
The~$\frac{n(n-1)}2$ random variables~$X_{i,j}=f(\gamma_i,\gamma_j)$ with~$i<j$ are identically distributed but not completely independent. In fact~$X_{i,j}$ is independent of~$X_{k,l}$ if and only if $\{i,j\}\cap\{k,l\}=\emptyset$. To apply the generalized weak law of large numbers, Theorem~\ref{thm:wlln_mostly_indep} below, we must check that the variables~$X_{i,j}$ are sufficiently independent. There are~$\frac{n(n-1)}{2} \sim n^2$ variables~$X_{i,j}$, which yield~$\sim n^4$ pairs~$(X_{i,j},X_{k,l})$, of which only~$\sim n^3$ are not independent. Therefore the proportion of nonindependent pairs $p\sim\frac{n^3}{n^4}\sim\frac 1n$ goes to zero as~$n\to \infty$. Thus, by Theorem~\ref{thm:wlln_mostly_indep}, the average value of the variables~$X_{i,j}$,
\[  \frac{\sum_{i<j}X_{i,j}}{\frac{n(n-1)}2}
= \frac {\sum_{i<j}\#(\gamma_i \cap \gamma_j)}{\frac{n(n-1)}2}
= \frac{\area(M,\W)}{\frac{n(n-1)}2} \]
converges in probability to the expected value, which, by the Santal\'o+Blaschke formula~\eqref{eq:croftonsantalo}, is equal to
\[ \Expe(X_{i,j})
= \iint_{\Gamma \times \Gamma}
 \#(\gamma \cap \gamma') \, \dProb(\gamma) \, \dProb(\gamma')
= \frac{2}{L^2} \,\pi\area(M).\]
\end{proof}

This concludes the proof of Theorem~\ref{theo:existence}.
\end{proof}

Let us prove the following generalization of the weak law of large numbers.

\begin{theorem}[Weak law of large numbers for identically distributed, mostly independent random variables]\label{thm:wlln_mostly_indep}
Fix a real valued random variable~$X$ with finite expected absolute value~$\Expe(|X|)<\infty$ and an integer~$n>0$.
Then the average $\overline X = \frac 1n \sum_i X_i$ of~$n$ random variables~$X_i$, each with the same distribution as~$X$, is near the expected value~$\Expe(X)$ with probability arbitrarily close to~$1$
if the proportion of nonindependent pairs
\[ p = \frac
{ \# \{(i,j)\mid X_i\text{ and }X_j\text{ are not independent} \} }
{n^2} \]
is small.
More precisely, for every $\eps,\delta>0$, there exists $p_0=p_0(X,\delta,\eps) >0$ such that if $p \leq p_0$, then
$\Prob \left( \left| \overline X - \Expe(X) \right| \geq \eps \right) \leq \delta$.
\end{theorem}

\begin{remark}
Note that we do not explicitly require~$n$ to be large, but this is generally necessary for~$p$ to be small, because each variable~$X_i$ is in general correlated with itself,\footnote{A random variable is independent of itself if and only if its probability distribution is concentrated in one value.} which implies that $p\geq\frac n{n^2}=\frac 1n$. If these are the only correlations and~$n$ goes to infinity, then $p=\frac 1n \to 0$ and therefore~$\overline X$ converges to~$\Expe(X)$ in probability. In this way, we recover the usual weak law of large numbers.
\end{remark}

\begin{proof} The proof is similar to the standard proof of the weak law of large numbers; see~\cite[Theorem~1.5.1]{tao2008strong} for instance.
It proceeds by cases; only the first one requires attention to the non-independent pairs.

{\it Case $\Expe(X^2)<\infty$ and $\Expe(X)=0$.} Fix $\eps>0$. We have to prove that the probability of deviation $\Prob \left(|\overline X|>\eps\right)$ gets arbitrarily low if~$p$ is sufficiently small.
To apply Chebyshev's inequality, we compute
\[ \Expe \left( \overline X^2 \right)
= \frac 1{n^2} \sum_i\sum_j \Expe(X_iX_j)
\leq p\, \Expe \left(X^2\right). \]
Here we used the Cauchy--Schwartz inequality $\Expe(X_iX_j) \leq \Expe(X^2)$ and the fact that $\Expe(X_iX_j)=\Expe(X_i)\Expe(X_j)=\Expe(X)^2=0$ if~$X_i$ and~$X_j$ are independent.
Applying Chebyshev's inequality, we obtain
\[ \Prob \left( |\overline X| \geq \eps \right)
  \leq \frac {\Expe \left( \overline X^2 \right) }{\eps^2}
  \leq \frac {p\, \Expe (X^2) }{\eps^2}
  \underset{p\to 0}{\longrightarrow} 0 \]
as we had to prove.

{\it Case $\Expe(X^2)<\infty$.} This case follows from the previous one applied to the random variable $Y=X-\Expe(X)$, which satisfies $\Expe(Y^2)<\infty$ and $\Expe(Y)=0$.

{\it General case $\Expe(|X|)<\infty$.} This case, which is not needed in this article, follows from a truncation argument as in the usual proof of the weak law of large numbers, given for instance in~\cite{tao2008strong}.

We proceed to the details. It is sufficient to show that
\begin{equation} \label{eq:wlln_target}
\Prob(|\overline X - \Expe(X)|\geq 3\eps) \leq 2\delta
\end{equation} if~$p$ is small enough with respect to~$\eps$,~$\delta$ and~$X$. We may assume $\delta \leq 1$.

We proceed as follows. For any cutoff value $M\geq 0$, we decompose the random variable~$X$ as a sum of a bounded part and a tail
\begin{equation}\label{eq:X-truncation}
X = X^{<M} + X^{\geq M}
\end{equation}
where the bounded part is
\[ X^{<M} = \mathbbm{1}_{|X|< M}\,X
= \begin{cases}
X &\text{if }|X| < M\\
0 &\text{otherwise},
\end{cases}\]
and the tail is
\[ X^{\geq M} = \mathbbm{1}_{|X|\geq M}\, X
= \begin{cases}
X &\text{if }|X| \geq M\\
0 &\text{otherwise}.
\end{cases}\]
In the same way we decompose the variables $X_i =  X_i^{<M} + X_i^{\geq M}$ and define two separate average values: one for the bounded parts, $\overline{X^{<M}} = \frac 1n \sum_i X_i^{<M}$, and one for the tails, $\overline{X^{\geq M}} = \frac 1n \sum_i X_i^{\geq M}$. These averages satisfy $\overline X = \overline{X^{<M}} + \overline{X^{\geq M}}$.

A key fact about the decomposition~\eqref{eq:X-truncation} is that the expected absolute value $\Expe\left(\left|X^{\geq M}\right|\right)$ of the tail part gets arbitrarily small if~$M$ is sufficiently large. This follows from the pointwise convergence $|X^{\geq M}|\to 0$ as $M\to+\infty$, which is dominated by~$|X|$, or from the formula
\[ \Expe(|X^{\geq M}|)
= \Expe(|X|^{\geq M})
= \int_M^{+\infty} x \diff\Prob_{|X|}(x), \]
where~$\Prob_{|X|}$ is the probability distribution of~$|X|$ on~$\R$.
We choose~$M$ so that
\begin{equation}\label{eq:tail_bound}
\Expe(|X^{\geq M}|) \leq \delta\eps.
\end{equation}
This implies that the average~$\overline{X^{\geq M}}$ of the tail parts also has small expected absolute value
\[ \Expe (| \overline{X^{\geq M}} |)
\leq \Expe(|X^{\geq M}|) \leq \delta\eps.\]
By Markov's inequality, this implies that~$\overline{X^{\geq M}}$ is small in absolute value with high probability
\begin{equation}\label{eq:tail_markov}
\Prob (| \overline{X^{\geq M}} | \geq \eps )
\leq \frac {\Expe (| \overline{ X^{\geq M}} |)} {\eps}
\leq \frac {\delta\eps} {\eps}
= \delta.
\end{equation}
Now, the bounded part~$X^{<M}$ has finite second moment
$\Expe((X^{<M})^2) < \infty$. Therefore, we may apply the previous case of the theorem, which yields
\begin{equation} \label{eq:bounded_markov}
\Prob\left( | \overline{X^{<M}} - \Expe (X^{<M}) | \geq \eps \right)
\leq \delta
\end{equation}
if~$p$ is small enough. Conveniently, $\Expe\left(X^{<M}\right)$ is near~$\Expe(X)$ because
\[ \left| \Expe(X) - \Expe\left(X^{<M}\right) \right|
= \left| \Expe \left( X^{\geq M}\right)\right|
\leq \Expe \left(\left|X^{\geq M}\right|\right)
\leq \delta\eps \leq \eps. \]
Here we used \eqref{eq:tail_bound} and the assumption $\delta\leq 1$.
Combining this with \eqref{eq:bounded_markov} and \eqref{eq:tail_markov} by the triangle inequality, the result \eqref{eq:wlln_target} follows.
\end{proof}

\section{Minimal area of disks: from discrete to Finsler metrics} \label{sec:discrete-to-smooth}

The goal of this section is to state a discrete version of the area lower bound on Finsler disks with minimizing interior geodesics and to show how to derive the area lower bound for Finsler metrics from its discrete version.

\medskip

Let us recall the area lower bound for Finsler metrics we want to prove.

\begin{theorem} \label{theo:6pi}
Let~$M$ be a self-reverse Finsler metric disk of radius~$r$ with minimizing interior geodesics.
Then the Holmes--Thompson area of~$M$ satisfies
\[   \area(M) \geq \frac{6}{\pi} \,r^2.   \]
\end{theorem}

In order to state the discrete version of this result, we need to introduce the notion of simple discrete metric disks.

\begin{definition} \label{def:sdmd}
A topological disk~$D$ with a quasi wall system~$\W$ is a \term{simple discrete metric disk} of radius~$r$ centered at an interior point $O \in D \setminus\W$ if the quasi wall system~$\W$ is simple (see Definition~\ref{def:wall}), all the points of~$D \setminus\W$ are at $d_\W$-distance at most~$r$ from~$O$ and all the points of $\partial D \setminus\W$ are at distance exactly~$r$ from~$O$.
\end{definition}

It is essential here to allow~$\W$ to be a quasi wall system rather than a wall system. Indeed, all points of $\W$ located on $\partial D$ necessarily have multiplicity~$2$. 

\medskip

The following result, which will be proved in the subsequent sections, can be seen as a discrete version of Theorem~\ref{theo:6pi}.

\begin{theorem} \label{theo:discrete-area}
The discrete area of every simple discrete metric disk~$(D,\W)$ of radius~$r$ satisfies
\[   \area(D,\W) \geq \frac{3}{2}\, r^2.   \]
Furthermore, the equality is attained.
\end{theorem}

Assuming this discrete area lower bound, we can derive Theorem~\ref{theo:6pi} as follows.

\begin{proof}[Proof of Theorem~\ref{theo:6pi} (assuming Theorem~\ref{theo:discrete-area})]
Let~$M$ be a Finsler disk of radius~$r$ centered at~$O$ with minimizing interior geodesics.
By Theorem~\ref{theo:existence}, for every~$\eps >0$, there exists a simple wall system~$\W_M$, made of~$n$ interior geodesics of~$M$, satisfying the estimates~\eqref{eq:distW} and~\eqref{eq:areaW}.
The simple wall system~$\W_M$ decomposes~$M$ into convex polygonal cells.
By definition, all the points in a cell are at the same distance from the center of~$M$ with respect to the discrete distance~$d_{\W_M}$.
Since~$M$ has minimizing interior geodesics, the geodesic rays of length~$r$ arising from its center~$O$ form a geodesic foliation~$\mathcal{F}$ of the punctured disk~$M \setminus \{ O \}$.
The sides of the cells of~$M$, which lie in the geodesics of~$\W_M$, are transverse to the foliation~$\mathcal{F}$, otherwise the origin~$O$ would lie in~$\W$.

\medskip

Consider a convex polygonal cell~$\Delta$ of~$M$ not containing~$O$.
Choose an arbitrary interior point of~$\Delta$ as its center.
Denote by~$d$ the $d_{\W_M}$-distance from~$O$ to the interior of the cell~$\Delta$.
The geodesic rays of the foliation intersecting~$\Delta$ form a spray~$\mathcal{F}_\Delta$, where each ray of~$\mathcal{F}_\Delta$ intersects~$\Delta$ along an interval with nonempty interior, except for the two extremal rays of the spray which intersect the convex polygonal cell~$\Delta$ at two vertices; see Figure~\ref{fig:spray}.
Denote by~$\beta_\Delta$ the broken line made of two segments joining the center of~$\Delta$ to these two extremal vertices.
Note that every geodesic ray of the spray~$\mathcal{F}_\Delta$ intersects the broken line~$\beta_\Delta$ at a single point; see Figure~\ref{fig:spray}.
Since the rays of the spray are length-minimizing with respect to~$d_{\W_M}$, see Proposition~\ref{prop:min}, all the cells intersecting the spray between~$O$ and~$\beta_\Delta$ are at $d_{\W_M}$-distance at most~$d$ from~$O$, and all the cells intersecting the spray after~$\beta_\Delta$ are at $d_{\W_M}$-distance at least~$d$ from~$O$.

\medskip

\begin{figure}[htbp!] 
\vspace{0.4cm}
\def\svgwidth{4.5cm}
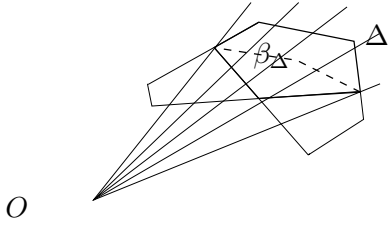 
\vspace{0.2cm}
\caption{Spray~$\mathcal{F}_\Delta$ intersecting the convex polygonal cell~$\Delta$} \label{fig:spray}
\end{figure}

Denote by~$r_0$ the integral part of~$n \left( \frac{2r}{L}-\eps \right)$.
By~\eqref{eq:distW}, every boundary point~$p \in \partial M \setminus\W_M$ is at $d_{\W_M}$-distance greater than~$r_0$ from~$O$, that is, $d_{\W_M}(O,p) > r_0$.
A cell of~$M$ whose interior points are at $d_{\W_M}$-distance~$r_0$ from~$O$ will be referred to as an outermost cell.
The broken lines~$\beta_\Delta$, where~$\Delta$ runs over all outermost cells of~$M$, form a piecewise geodesic closed curve delimiting a topological disk~$D \subseteq M$ containing~$O$. This curve can be smoothed to ensure that~$D$ is a smoothly bounded manifold.
The restriction $\W= \W_M \cap D$ of~$\W_M$ to~$D$ defines a simple quasi wall system on~$D$.
By construction, all the points of~$D \setminus\W$ are at $d_\W$-distance at most~$r_0$ from~$O$ and all the points of $\partial D \setminus\W$ are at distance exactly~$r_0$ from~$O$.
Hence,~$(D,\W)$ is a simple discrete metric disk of radius~$r_0$.
By Theorem~\ref{theo:discrete-area} and by definition of the discrete area~\eqref{eq:area}, we have
\[
\area(M,\W_M) \geq \area(D,\W) \geq \frac{3}{2} \, r_0^2.
\]
Dividing by~$n^2$, using~\eqref{eq:distW} and~\eqref{eq:areaW}, and letting~$\eps$ go to zero, we obtain
\[
\frac{\pi}{L^2} \, \area(M) \geq \frac{3}{2} \left( \frac{2}{L} r \right)^2.
\]
Hence, $\area(M) \geq \frac{6}{\pi} \, r^2$.
\end{proof}

Sections~\ref{sec:semi-circle}--\ref{sec:proof} are devoted to the proof of Theorem~\ref{theo:discrete-area}.

\section{Quasi wall systems and interval families} \label{sec:semi-circle}

In this section, we show how to encode a simple discrete disk as a 1-dimensional object.

\medskip

We start by proving the following basic fact about simple discrete metrics.

\begin{proposition}\label{thm:distance_simple}
Let~$D$ be a disk with a simple quasi wall system~$\W$. Then
\begin{equation} \label{eq:distance_simple}
d_\W(x,y) = \text{number of walls of } \W
\text{ that separate } x \text{ from } y.
\end{equation}
for any two points~$x,y\in D\setminus\W$.
\end{proposition}

Note that if~$D$ is a Finsler disk with minimizing interior geodesics and~$\W$ is geodesic, then this proposition follows from Proposition~\ref{prop:min}.

\begin{proof} It is clear that
\[ d_\W(x,y) \geq
\text{ number of walls of } \W \text{ that separate } x \text{ from } y. \]
To prove the reverse inequality we will show the following.
\begin{claim} There exists a smooth path~$\alpha$ from~$x$ to~$y$ that is in general position with respect to $\W\cup\partial D$ and crosses each wall of~$\W$ at most once.
\end{claim}
Here, we say that a smooth curve~$\alpha$ is in general position with respect to an immersed $1$-submanifold~$N$ if it is regular, transverse to~$N$ and avoids the self-intersections of~$N$. If~$\alpha$ is piecewise smooth, we require in addition that none of its non-smooth points lie in~$N$.

The claim is a version of Levi's extension (or enlargement) lemma for pseudoline arrangements. This version concerns arrangements on a disk, rather than on the projective plane as in the more standard version of the lemma (found e.g. in ~\cite[Thm. 5.1.1]{felsner2018pseudoline}).

We prove the claim by induction on the number of walls. Suppose the claim is valid for any quasi wall system~$\W$ made of~$n$ walls.
Consider a simple quasi wall system~$\W'$ obtained by adding an extra~$w'$ to~$\W$. By inductive hypothesis, there is a smooth path~$\alpha$ that satisfies all the conditions of the claim with respect to~$\W$. By perturbing~$\alpha$, we ensure that it is transverse to~$w'$ as well. If~$\alpha$ crosses~$w'$ at most once, then we are done. Otherwise, let~$x'$ and~$y'$ be the first and last points of~$\alpha$ where~$\alpha$ crosses~$w'$. Note that they are generic points of~$w'$: they are neither on~$\W$, nor on~$\partial D$. Replace the segment of~$\alpha$ from~$x'$ to~$y'$ by the segment $[x',y']$ of~$w'$, and let~$\alpha'$ be the resulting curve. We claim that~$\alpha'$ is a piecewise smooth curve, in general position with respect to~$\W$, that crosses each wall of~$\W$ at most once. This is because the segment $[x',y']$ that we inserted only crosses the walls of~$\W$ that separate~$x'$ from~$y'$ (since it is part of a wall of the simple quasi wall system~$\W'$), and these walls are necessarily crossed as well by the piece of~$\alpha$ between~$x'$ and~$y'$ that we replaced.

The next step is to perturb the curve~$\alpha'$ so that the segment $[x',y']$ is displaced sideways and away from~$w'$ and the resulting curve~$\alpha''$ is in general position with respect to~$\W\cup\partial D$ and crosses~$\W$ the same number of times as~$\alpha'$ does, and, in addition, is transverse to~$w'$ and crosses~$w'$ at most once. Thus,~$\alpha''$ is in general position with respect to $\W'\cup\partial D$ and crosses each wall of~$\W'$ at most once, but is non-smooth at two points. To make it smooth, we modify it near these two points.
\end{proof}

Let~$(D,\W)$ be a simple discrete disk of radius~$r$ and center~$O$, see Definition~\ref{def:sdmd}.
Identify the boundary~$\partial D$ with the circle~$S^1$, and 
identify the punctured disk~$D \setminus \{O\}$ with the flat cylinder~$\C=S^1 \times [0,\infty)$. Under this identification, the point~$O$ of~$D$ corresponds to the point at infinity in the one-point compactification of the cylinder~$\C$. Note that the universal cover of~$\C$ is the half plane $\HH = \R\times[0,+\infty)$.

\begin{definition} \label{def:cover}
Given a simple arc~$\alpha$ in the cylinder~$M = D \setminus \{ O \}$ (or in the half plane $M=\R \times [0,+\infty)$) with endpoints on the boundary~$\partial M$, denote by~$\overline{\alpha}$ the segment of~$\partial D$ with the same endpoints, homotopic to~$\alpha$ in~$M$.
The arc~$\alpha$ \term{covers a point}~$p$ of~$\partial M$ if~$p$ lies in~$\overline{\alpha}$.
Similarly, the arc~$\alpha$ \term{covers another arc}~$\beta$ if~$\overline{\beta}$ lies in~$\overline{\alpha}$.
Two arcs~$\alpha$ and~$\alpha'$ are \term{adjacent} if the intervals~$\overline{\alpha}$ and~$\overline{\alpha'}$ are adjacent, meaning that they have exactly one point in common.
\end{definition}

\begin{definition}\label{def:standard}
An arc in the flat cylinder $M = S^1\times[0,\infty)$ (or in the half plane plane $M = \R \times [0,+\infty)$) with endpoints on the boundary~$\partial M$ is \term{standard} if it consists of a segment of slope~$1$ followed by a segment of slope~$-1$; see Figure~\ref{fig:2inter}.
A quasi wall system~$\W$ is standard if its walls are standard arcs.
For two boundary points $a,b\in S^1=\partial M$ (or $a<b\in\R$ if~$M$ is the half plane), we denote by~$[a,b]$ the arc of~$S^1$ that goes from~$a$ to~$b$ in the positive (\ie, counterclockwise) sense,
and we denote by~$\widehat{ab}$ the standard arc in~$M$ that is homotopic to~$[a,b]$.
\end{definition}

Let~$(D,\W)$ be a simple discrete disk of radius~$r$ centered at~$O$. Denote by~$\II=\II_\W$ the set of boundary intervals~$\overline\alpha$ homotopic to the walls~$\alpha$ of~$\W$.
The family~$\II$ of intervals of~$S^1$ contains all the information about~$\W$ that is relevant to our problem of finding simple discrete disks of minimum area. For instance, two walls~$\alpha,\beta$ of~$\W$ meet on~$\partial D$ if and only if the intervals~$\overline\alpha$,~$\overline\beta$ have a common endpoint. That is,
\begin{equation}\label{eq:intervals_meet}
\# I_{\partial D}(\alpha,\beta)=1 \iff
\# (\partial\overline\alpha \cap \partial\overline\beta) = 1.
\end{equation}
Furthermore, assuming~$\overline\alpha$ and~$\overline\beta$ have no common endpoints, the arcs~$\alpha$ and~$\beta$ cross in the interior of~$D$ if and only if the interval~$\overline\alpha$ contains \emph{exactly one} endpoint of~$\overline\beta$. That is,
\begin{equation}\label{eq:intervals_cross}
\# I_{\Int D}(\alpha,\beta)=1 \iff
\# (\overline\alpha \cap \partial\overline\beta) = 1.
\end{equation}
One consequence of these formulas is that the discrete area of~$(D,\W)$ given by~\eqref{eq:area} may be computed from~$\II$.

\medskip

The following result characterizes the relation between the quasi wall system~$\W$ and the interval family~$\II$.
Before stating this result, we need to introduce a definition.
A point~$p$ of~$S^1$ is \term{generic} with respect to a finite interval family~$\II$ of~$S^1$ if~$p$ is not an endpoint of any interval of~$\II$.
Alternatively, the endpoints of the intervals of~$\II$ are the non-generic points of~$S^1$.

\begin{proposition}\label{thm:simple_interval_system}
Let~$(D,W)$ be a simple discrete disk of radius~$r$ centered at~$O$. The family~$\II = \II_\W$ of intervals of~$S^1$ has the following properties:
\begin{enumerate}
\item\label{item:no_two_cover} no pair of intervals of~$\II$ cover~$S^1$;
\item\label{item:r_cover} every generic point of~$S^1$ is contained in exactly~$r$ intervals of~$\II$;
\item\label{item:interval_genericity} every non-generic point of~$S^1$ is an endpoint of exactly two, adjacent intervals of~$\II$.
\end{enumerate}
Moreover, if a finite family~$\II$ of intervals of~$S^1$ satisfies the conditions~\eqref{item:no_two_cover}--\eqref{item:interval_genericity}, then~$\II=\II_\W$ for some quasi wall system~$\W$ that makes~$D$ a simple discrete metric disk of radius~$r$ and center~$O$. For instance, one may let~$\W$ be the unique standard quasi wall system homotopic to~$\II$ on~$D\setminus\{O\}$.
\end{proposition}

\begin{proof}\mbox{}\\
\noindent\eqref{item:no_two_cover} If two intervals $\overline\alpha,\overline\beta \in \II$ cover~$S^1$, then the corresponding walls~$\alpha,\,\beta$ of~$\W$ would form a bigon containing the point~$O$, which implies they cross twice, contradicting the hypothesis that~$\W$ is simple.\\
\noindent\eqref{item:r_cover} Consider a generic point~$p\in S^1$.
Since~$\W$ is a simple quasi wall system on~$D$, the distance between any pair of points of~$D$ is the number of walls that separate them; see Proposition~\ref{thm:distance_simple}. On the other hand, the walls that separate~$O$ from~$p$ are the walls that cover~$p$.
Hence the result.\\
\noindent\eqref{item:interval_genericity} This follows from the previous property: if~$p\in S^1$ is the endpoint of some interval~$\overline\alpha\in \II$, it must also be the startpoint of some other interval so that every generic point near~$p$ is contained in the same number~$r$ of intervals of~$\II$. This means that~$p$ is the endpoint of two walls, and it cannot be the endpoint of more walls because~$\W$ can only have simple self-intersections on~$\partial D$ since it is a quasi wall system; see Definition~\ref{def:wall}.

Now, let~$\II$ be a finite family of intervals of~$S^1$ satisfying conditions~\eqref{item:no_two_cover}--\eqref{item:interval_genericity}, and let~$\W$ be the unique standard quasi wall system homotopic to~$\II$ on~$D\setminus\{O\}$. Clearly,~$\W$ is a quasi wall system, and it is simple because it is made of arcs that intersect each other at most once. Also, every point $p\in D\setminus\W$ is at distance at most~$r$ from~$O$, and exactly~$r$ if~$p\in\partial D$. (A shortest path is the vertical ray from~$p$ to~$O$.) This shows that~$(D,\W)$ is a simple discrete disk of radius~$r$ centered at~$O$.
\end{proof}

\section{Inadmissible configurations in a minimal simple disk} \label{sec:inadmissible}

In this section, we rule out some intersection patterns for an extremal quasi wall system on a disk. 

\medskip

Consider a quasi wall system~$\W$ on~$D$ defining a simple discrete metric disk of radius~$r$ \emph{with minimal discrete area}.
By Proposition~\eqref{thm:simple_interval_system}, we can assume that~$\W$ is formed of standard arcs; see Definition~\ref{def:standard}.

\begin{lemma} \label{lem:nocover}
No arc of~$\W$ covers two (possibly adjacent) intersecting arcs of~$\W$.
\end{lemma}

\begin{proof}
By contradiction, suppose that an arc~$\gamma$ of~$\W$ covers two intersecting arcs~$\alpha=\widehat{ac}$ and~$\beta=\widehat{bd}$ of~$\W$. Switching the roles of the two arcs if necessary, we may assume that the points~$a,\,b,\,c,\,d$ appear in that order in the interval~$\overline\gamma$ (with possibly~$b=c$).
See Figure~\ref{fig:2inter}.
Let~$\W'$ be the collection of curves obtained from~$\W$ by replacing~$\alpha$ and~$\beta$ with the standard arcs~$\alpha'=\widehat{ad}$ and~$\beta'=\widehat{bc}$ (with no~$\beta'$ if~$b=c$).
See Figure~\ref{fig:2inter}.
Note that, like~$\W$, the immersed $1$-submanifold~$\W'$ is a quasi wall system on~$D$. Moreover, we claim that~$\W'$ also makes~$D$ a simple discrete metric disk of radius~$r$ centered at~$O$. This is because none of the properties~\eqref{item:no_two_cover}--\eqref{item:interval_genericity} of Proposition~\ref{thm:simple_interval_system} is affected by the replacement. For instance, there is no arc~$\delta$ of~$\W'$ such that the intervals~$\overline\delta$ and~$\overline{\alpha'}$ cover the boundary~$\partial D$, because in that case~$\overline\delta$ and~$\overline\gamma$ would also cover~$\partial D$, however the arcs~$\overline\delta$ and~$\overline\gamma$ are already present in~$\W$, contradicting by Proposition~\ref{thm:simple_interval_system} the fact that~$\W$ is simple. Also, the fact that every generic point of~$\partial D$ is covered by exactly~$r$ arcs of the quasi wall system is clearly maintained, as well as the fact that each non-generic boundary point is the common endpoint of two adjacent walls.

\begin{figure}[htbp!] 
\vspace{.2cm}
\def\svgwidth{8cm}
\begingroup%
  \makeatletter%
  \providecommand\color[2][]{%
    \errmessage{(Inkscape) Color is used for the text in Inkscape, but the package 'color.sty' is not loaded}%
    \renewcommand\color[2][]{}%
  }%
  \providecommand\transparent[1]{%
    \errmessage{(Inkscape) Transparency is used (non-zero) for the text in Inkscape, but the package 'transparent.sty' is not loaded}%
    \renewcommand\transparent[1]{}%
  }%
  \providecommand\rotatebox[2]{#2}%
  \newcommand*\fsize{\dimexpr\f@size pt\relax}%
  \newcommand*\lineheight[1]{\fontsize{\fsize}{#1\fsize}\selectfont}%
  \ifx\svgwidth\undefined%
    \setlength{\unitlength}{226.77165354bp}%
    \ifx\svgscale\undefined%
      \relax%
    \else%
      \setlength{\unitlength}{\unitlength * \real{\svgscale}}%
    \fi%
  \else%
    \setlength{\unitlength}{\svgwidth}%
  \fi%
  \global\let\svgwidth\undefined%
  \global\let\svgscale\undefined%
  \makeatother%
  \begin{picture}(1,0.5)%
    \lineheight{1}%
    \setlength\tabcolsep{0pt}%
    \put(0,0){\includegraphics[width=\unitlength,page=1]{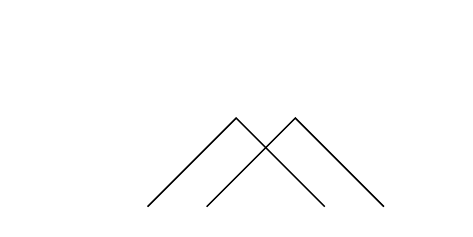}}%
    \put(0.48635974,0.26817467){\makebox(0,0)[lt]{\lineheight{1.25}\smash{\begin{tabular}[t]{l}$\alpha$\end{tabular}}}}%
    \put(0.60969367,0.26814771){\makebox(0,0)[lt]{\lineheight{1.25}\smash{\begin{tabular}[t]{l}$\beta$\end{tabular}}}}%
    \put(0,0){\includegraphics[width=\unitlength,page=2]{fig_2inter_old.pdf}}%
    \put(0.54662333,0.45949976){\color[rgb]{0,0,0}\makebox(0,0)[lt]{\lineheight{1.25}\smash{\begin{tabular}[t]{l}$\gamma$\end{tabular}}}}%
    \put(0.28849266,0.3253405){\makebox(0,0)[lt]{\lineheight{1.25}\smash{\begin{tabular}[t]{l}$\delta$\end{tabular}}}}%
    \put(0.425,0.0125){\makebox(0,0)[lt]{\lineheight{1.25}\smash{\begin{tabular}[t]{l}$b$\end{tabular}}}}%
    \put(0.675,0.0125){\makebox(0,0)[lt]{\lineheight{1.25}\smash{\begin{tabular}[t]{l}$c$\end{tabular}}}}%
    \put(0.8,0.0125){\makebox(0,0)[lt]{\lineheight{1.25}\smash{\begin{tabular}[t]{l}$d$\end{tabular}}}}%
    \put(0.3,0.0125){\makebox(0,0)[lt]{\lineheight{1.25}\smash{\begin{tabular}[t]{l}$a$\end{tabular}}}}%
  \end{picture}%
\endgroup%

\def\svgwidth{8cm}
\begingroup%
  \makeatletter%
  \providecommand\color[2][]{%
    \errmessage{(Inkscape) Color is used for the text in Inkscape, but the package 'color.sty' is not loaded}%
    \renewcommand\color[2][]{}%
  }%
  \providecommand\transparent[1]{%
    \errmessage{(Inkscape) Transparency is used (non-zero) for the text in Inkscape, but the package 'transparent.sty' is not loaded}%
    \renewcommand\transparent[1]{}%
  }%
  \providecommand\rotatebox[2]{#2}%
  \newcommand*\fsize{\dimexpr\f@size pt\relax}%
  \newcommand*\lineheight[1]{\fontsize{\fsize}{#1\fsize}\selectfont}%
  \ifx\svgwidth\undefined%
    \setlength{\unitlength}{226.77165354bp}%
    \ifx\svgscale\undefined%
      \relax%
    \else%
      \setlength{\unitlength}{\unitlength * \real{\svgscale}}%
    \fi%
  \else%
    \setlength{\unitlength}{\svgwidth}%
  \fi%
  \global\let\svgwidth\undefined%
  \global\let\svgscale\undefined%
  \makeatother%
  \begin{picture}(1,0.5)%
    \lineheight{1}%
    \setlength\tabcolsep{0pt}%
    \put(0,0){\includegraphics[width=\unitlength,page=1]{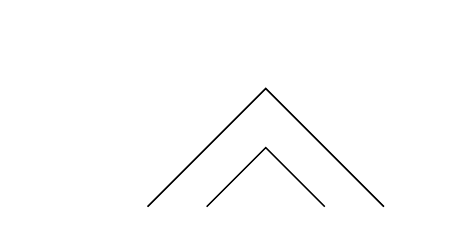}}%
    \put(0.54322147,0.33201536){\makebox(0,0)[lt]{\lineheight{1.25}\smash{\begin{tabular}[t]{l}$\alpha'$\end{tabular}}}}%
    \put(0.54442086,0.20641048){\makebox(0,0)[lt]{\lineheight{1.25}\smash{\begin{tabular}[t]{l}$\beta'$\end{tabular}}}}%
    \put(0,0){\includegraphics[width=\unitlength,page=2]{fig_2inter_new.pdf}}%
    \put(0.54662333,0.45949976){\color[rgb]{0,0,0}\makebox(0,0)[lt]{\lineheight{1.25}\smash{\begin{tabular}[t]{l}$\gamma$\end{tabular}}}}%
    \put(0.28849266,0.3253405){\makebox(0,0)[lt]{\lineheight{1.25}\smash{\begin{tabular}[t]{l}$\delta$\end{tabular}}}}%
    \put(0.425,0.0125){\makebox(0,0)[lt]{\lineheight{1.25}\smash{\begin{tabular}[t]{l}$b$\end{tabular}}}}%
    \put(0.675,0.0125){\makebox(0,0)[lt]{\lineheight{1.25}\smash{\begin{tabular}[t]{l}$c$\end{tabular}}}}%
    \put(0.8,0.0125){\makebox(0,0)[lt]{\lineheight{1.25}\smash{\begin{tabular}[t]{l}$d$\end{tabular}}}}%
    \put(0.3,0.0125){\makebox(0,0)[lt]{\lineheight{1.25}\smash{\begin{tabular}[t]{l}$a$\end{tabular}}}}%
  \end{picture}%
\endgroup%
 
\caption{Replacing two intersecting arcs covered by a third arc} \label{fig:2inter}
\end{figure}

\medskip

Let us show that the area of~$(D,\W')$ is less than the area of~$(D,\W)$ by comparing the number of self intersections of the quasi wall systems~$\W$ and~$\W'$ according to the discrete area formula~\eqref{eq:area}.
First, observe that every pair of arcs of~$\W$ different from~$\alpha$ and~$\beta$ belongs to~$\W'$.
Therefore, these pairs of arcs give the same contribution to the discrete areas of~$\W$ and~$\W'$.
Let~$\delta=\widehat{pq}$ be an arc of~$\W$ different from~$\alpha$ and~$\beta$.
By considering cases regarding the location of the endpoints~$p$ and~$q$ with respect to the points~$a$,~$b$,~$c$ and~$d$, we see that 
\[     \#I_{\Int D + \frac 12 \partial D}(\delta,\alpha'\cup\beta')
  \leq \#I_{\Int D + \frac 12 \partial D}(\delta,\alpha\cup\beta). \]
In fact, equality holds unless~$p$ and~$q$ lie in the interiors of~$[a,b]$ and~$[c,d]$, in which case the inequality is strict.
Finally, note that
\[ \#I_{\Int D}(\alpha',\beta') = 0 \quad\text{and}\quad
   \#I_{\Int D}(\alpha,\beta) = 1.\]
We conclude that
\[ \area (D,\W') \leq \area (D,\W) - 1, \]
which contradicts the minimiality of the discrete area of~$(D,\W)$.
\end{proof}

\begin{lemma} \label{lem:nointer}
No arc of~$\W$ intersects two adjacent arcs of~$\W$.
\end{lemma}

\begin{proof}
By contradiction, suppose that an arc~$\gamma$ of~$\W$ intersects two adjacent arcs~$\alpha$ and~$\beta$ of~$\W$. 
We choose~$\gamma$ so that it is minimal with respect to the covering relation, among arcs that intersects~$\alpha$ and~$\beta$ (\ie, no arc of~$\W$ covered by~$\gamma$ intersects~$\alpha$ and~$\beta$).
Denote by~$a,b,c,d,e$ the endpoints of the three arcs, in the order in which they are found on the interval $\overline\alpha\cup\overline\beta$. Thus,~$\alpha=\widehat{ac}$, $\beta=\widehat{ce}$ and~$\gamma = \widehat{bd}$, and no arc of~$\W$ that covers~$c$ is covered by~$\gamma$ (other than~$\gamma$ itself).
See Figure~\ref{fig:nointer}.
Let~$c^-$ and~$c^+$ be two points of~$\partial D$ close to~$c$ such that~$[c^-,c^+]\cap\partial\W = \{c\}$.
Let~$\W'$ be the collection of curves obtained from~$\W$ by replacing the three arcs~$\alpha$,~$\beta$ and~$\gamma$ with the four arcs~
$\alpha'=\widehat{ac^+}$,~
$\beta'=\widehat{c^-e}$,~
$\gamma^-=\widehat{bc^-}$ and~
$\gamma^+=\widehat{c^+d}$.
See Figure~\ref{fig:nointer}.

Note that~$\W'$ is a quasi wall system on the disk~$D$. In fact,~$\W'$ makes~$D$ a simple discrete disk of radius~$r$ centered at~$O$. To see this we argue as in the proof of Lemma~\ref{lem:nocover}. By Proposition~\ref{thm:simple_interval_system}, it is enough to check that the family~$\II=\II_{\W'}$ of boundary segments~$\overline\delta$ corresponding to the walls~$\delta$ of~$\W'$ satisfies the properties~\eqref{item:no_two_cover}--\eqref{item:interval_genericity} of Proposition~\ref{thm:simple_interval_system}.
To check Property~\eqref{item:r_cover} (that each generic point of~$\partial D$ is covered~$r$ times by the walls of~$\W'$) note that both~$\alpha\cup\beta\cup\gamma$ and~$\alpha'\cup\beta'\cup\gamma^-\cup\gamma^+$ cover twice the generic points of~$[b,d]$ and once the remaining generic points of~$[a,e]$.
Property~\eqref{item:interval_genericity} regarding non-generic boundary points is also maintained, with the wall endpoint~$c$ replaced by the two points~$c^-$ and~$c^+$.
Finally, to check the property~\eqref{item:no_two_cover}, suppose~$\delta$ and~$\eps$ are two arcs of~$\W'$ that cover the whole boundary~$\partial D$. It is impossible that both~$\delta$ and~$\eps$ are among the new arcs~$\alpha'$, $\beta'$ and~$\gamma^\pm$ because that would mean that~$\alpha$ and~$\beta$ already cover~$\partial D$, contradicting the fact that~$\W$ is simple. Similarly, the arcs~$\delta$ and~$\eps$ cannot be both among the unchanged arcs (those in~$\W\cap\W'$) either, otherwise~$\W$ would not be simple.
Therefore,~$\delta$ is one of the unchanged arcs and~$\eps$ is one of the new arcs~$\alpha'$, $\beta'$,~$\gamma^\pm$. In the case $\eps=\alpha'$, we see that~$\delta$ and~$\alpha'$ cannot cover~$\partial D$ since this would imply that~$\delta$ and~$\alpha$ already cover~$\partial D$. This is because $\overline{\alpha'}\setminus\overline\alpha$ is contained in the interval~$[c^-,c^+]$ which contains no endpoints of~$\delta$ since $[c^-,c^+]\cap\W = \{c\}$.
The case $\eps = \beta'$ is analogous and the cases $\eps = \gamma^\pm$ are easier to rule out since the arcs~$\gamma^\pm$ are covered by~$\gamma$. We conclude that the property~\eqref{item:no_two_cover} is satisfied, thus~$(D,\W')$ is a simple discrete metric disk of radius~$r$.

\begin{figure}[htbp] 
\def\svgwidth{8cm}
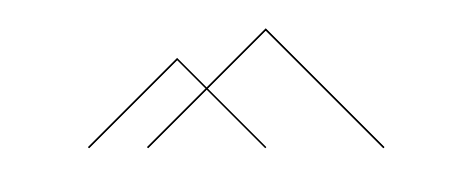
\def\svgwidth{8cm}
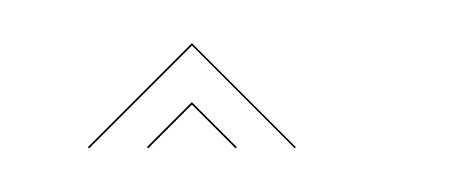
\caption{Replacing a configuration of one arc intersecting two adjacent arcs.}
\label{fig:nointer}
\end{figure}

\medskip

Let us show that the area of~$(D,\W')$ is less than the area of~$(D,\W)$. Again, we use the discrete area formula~\eqref{eq:area}, which says
\[ \area(D,\W) = \sum_{ \{\delta,\eps\} }
\# I_{\Int D + \frac 12 \partial D} (\delta,\eps) \]
where the sum runs over pairs $\{\delta,\eps\}$ of different walls of~$\W$.
The pairs $\{\delta,\eps\}$ of walls that are contained in $\W\cap\W'$ make the same contribution to~$\area(D,\W)$ and to~$\area(D,\W')$. To evaluate the contribution of pairs $\{\delta,\eps\}$ with $\delta\in W\cap\W'$ and $\eps\not\in\W\cap\W'$, we note that any arc $\delta=\widehat{pq}$ with no endpoints in~$[c^-,c^+]$ satisfies
\[ \#I_{\Int D + \frac 12 \partial D}
(\delta,\alpha' \cup \beta' \cup \gamma^+ \cup \gamma^-) = 
\#I_{\Int D + \frac 12 \partial D}
(\delta,\alpha\cup\beta\cup\gamma) \] 
unless~$p\in[b,c]$ and~$q\in[c,d]$.
This is seen by considering case by case the possible locations of~$p$ and~$q$ with respect to $a,\,b,\,c,\,d\,e$.
The equality holds for all arcs $\delta=\widehat{pq}\in\W\cap\W'$, because the exceptional case $p\in[b,c]$ and $q\in[c,d]$ is excluded by how~$\gamma$ was chosen: the arc $\gamma = \widehat{cd}$ covers no other arc $\delta=\widehat{pq}$ of~$\W$ that in turn covers~$c$. Finally, to compute the contribution of the pairs $\{\delta,\eps\}$ where none of the two arcs~$\delta$ and~$\eps$ is in~$\W\cap\W'$, we note that
\[ \# I_{\Int D + \frac 12 \partial D}
(\alpha' \cup \beta' \cup \gamma^- \cup \gamma^+) = 2 \]
while
\[ \# I_{\Int D + \frac 12 \partial D}
(\alpha \cup \beta \cup \gamma) = \frac 52.\]
We conclude that $\area(D,\W')=\area(D,\W)-\frac 12$,
contradicting the minimality of~$\W$.
\end{proof}

\section{Pairs of adjacent arcs} \label{sec:adjacent}

In this section we show that the sequences of adjacent arcs in an extremal quasi wall system on a disk have a periodic structure.

\medskip

Consider a quasi wall system~$\W$ on the disk~$D$, made of standard arcs, defining a simple discrete metric disk of radius~$r$ centered at~$O$ with minimal discrete area as in Section~\ref{sec:inadmissible}. 
Recall that the upper half plane~$\HH = \R \times [0,\infty)$ is the universal cover of the cylinder~$\C=S^1 \times [0,\infty) = D\setminus \{O\}$. We identify its boundary~$\partial\HH$ with the real line~$\R$.
Let~${\W_\HH}$ be the quasi wall system on~$\HH$ formed of all the lifts of the arcs of~$\W$. 

Since~$D$ is a disk of radius~$r$, it follows that every generic point of~$\partial \HH$ is covered by exactly~$r$ arcs of~$\W_\HH$.
To ensure this uniform coverage, each endpoint of an arc must be the startpoint of another arc, and thus each arc of~$\W_\HH$ belongs to a bi-infinite sequence of consecutive arcs, called a ``strand''~of $\W_\HH$.

\begin{definition} \label{def:strand}
A \term{strand} of $\W_\HH$ is a bi-infinite sequence~$(\alpha_i)_{i\in\Z}$ of consecutive arcs of~$\W_\HH$ of the form
\[ \alpha_i = \widehat{a_ia_{i+1}}. \]
The points~$a_i$ where the strand~$(\alpha_i)_i$ meets the boundary~$\partial\HH$ are called the \term{stops} of the strand.
The \term{width} of an arc~$\alpha_i$ is the number $a_{i+1}-a_i$.
\end{definition}

Since each strand of arcs covers the generic points of $\partial\HH$ once, it follows that the quasi wall system ~$\W_\HH$ is composed of exactly~$r$ strands.

\medskip

The following result describes how each strand intersects a pair of adjacent arcs of~$\W_\HH$.

\begin{lemma} \label{lem:zerotwo}
Let~$\alpha_0=\widehat{a_0a_1}$ and~$\alpha_1=\widehat{a_1a_2}$ be two adjacent arcs of~${\W_\HH}$.
Then every strand of~${\W_\HH}$ has exactly one arc with endpoints on the boundary interval $I=[a_0,a_2)$. This arc is covered by~$\alpha_0$ or by~$\alpha_1$.
\end{lemma}

\begin{proof} The strand that contains the arcs~$\alpha_0$ and~$\alpha_1$ clearly satisfies the proposition. Thus let $(\beta_i)_{i\in\Z}$ be any other strand of~${\W_\HH}$, numbered so that the arc~$\beta_0$ covers the point~$a_1$.
This strand has a stop in~$I$, otherwise~$\beta_0$ would cover the two adjacent arcs~$\alpha_0$ and~$\alpha_1$, in contradiction with Lemma~\ref{lem:nocover}. Also, the strand~$(\beta_i)_i$ cannot have stops in both intervals $[a_0,a_1)$ and $[a_1,a_2)$, otherwise the arc~$\beta_0$ would intersect the two adjacent arcs~$\alpha_0$ and~$\alpha_1$, in contradiction with Lemma~\ref{lem:nointer}. Thus the strand~$(\beta_i)_i$ has stops in exactly one of the intervals $[a_0,a_1)$ and $[a_1,a_2)$, say, the second one; see Figure~\ref{fig:leaping}. Furthermore, it cannot have just one stop in this interval, otherwise the two adjacent arcs $\beta_0,\,\beta_1$ that share this stop would intersect~$\alpha_1$, in contradiction with Lemma~\ref{lem:nointer}. Also, it cannot have three stops in the interval, otherwise the adjacent arcs~$\beta_1$ and~$\beta_2$ would be covered by~$\alpha_1$, in contradiction with Lemma~\ref{lem:nocover}. We conclude that the strand~$(\beta_i)_i$ has exactly two stops (and therefore one arc) in the interval $[a_0,a_2)$, and both of these stops are covered by one of the arcs~$\alpha_0$ or~$\alpha_1$.
See Figure~\ref{fig:leaping}.

\begin{figure}[htbp]
\vspace{-2mm}
\def\svgwidth{128mm}
\begingroup%
  \makeatletter%
  \providecommand\color[2][]{%
    \errmessage{(Inkscape) Color is used for the text in Inkscape, but the package 'color.sty' is not loaded}%
    \renewcommand\color[2][]{}%
  }%
  \providecommand\transparent[1]{%
    \errmessage{(Inkscape) Transparency is used (non-zero) for the text in Inkscape, but the package 'transparent.sty' is not loaded}%
    \renewcommand\transparent[1]{}%
  }%
  \providecommand\rotatebox[2]{#2}%
  \newcommand*\fsize{\dimexpr\f@size pt\relax}%
  \newcommand*\lineheight[1]{\fontsize{\fsize}{#1\fsize}\selectfont}%
  \ifx\svgwidth\undefined%
    \setlength{\unitlength}{532.91338583bp}%
    \ifx\svgscale\undefined%
      \relax%
    \else%
      \setlength{\unitlength}{\unitlength * \real{\svgscale}}%
    \fi%
  \else%
    \setlength{\unitlength}{\svgwidth}%
  \fi%
  \global\let\svgwidth\undefined%
  \global\let\svgscale\undefined%
  \makeatother%
  \begin{picture}(1,0.21276596)%
    \lineheight{1}%
    \setlength\tabcolsep{0pt}%
    \put(0,0){\includegraphics[width=\unitlength,page=1]{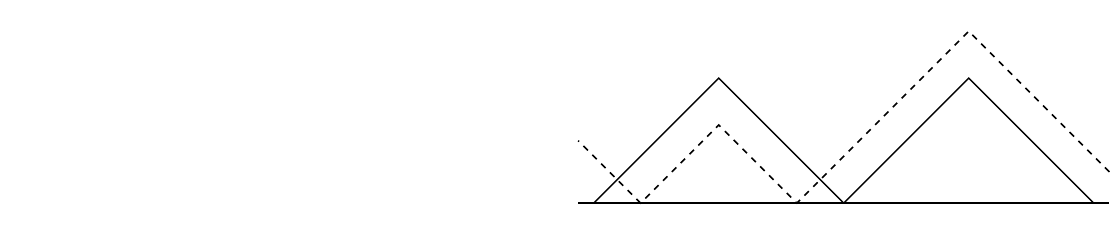}}%
    \put(0.07527308,0.06988563){\makebox(0,0)[lt]{\lineheight{1.25}\smash{\begin{tabular}[t]{l}$\alpha_0$\end{tabular}}}}%
    \put(0.27732952,0.11411688){\makebox(0,0)[lt]{\lineheight{1.25}\smash{\begin{tabular}[t]{l}$\alpha_1$\end{tabular}}}}%
    \put(0.03426683,0.13521578){\makebox(0,0)[lt]{\lineheight{1.25}\smash{\begin{tabular}[t]{l}$\beta_0$\end{tabular}}}}%
    \put(0,0){\includegraphics[width=\unitlength,page=2]{fig_leaping.pdf}}%
    \put(-0.00437719,0.00531915){\makebox(0,0)[lt]{\lineheight{1.25}\smash{\begin{tabular}[t]{l}$a_0$\end{tabular}}}}%
    \put(0.22884212,0.00531915){\makebox(0,0)[lt]{\lineheight{1.25}\smash{\begin{tabular}[t]{l}$a_1$\end{tabular}}}}%
    \put(0.45804049,0.0063244){\makebox(0,0)[lt]{\lineheight{1.25}\smash{\begin{tabular}[t]{l}$a_2$\end{tabular}}}}%
    \put(0.33182255,0.05177954){\makebox(0,0)[lt]{\lineheight{1.25}\smash{\begin{tabular}[t]{l}$\beta_1$\end{tabular}}}}%
  \end{picture}%
\endgroup%

\caption{Leaping over every other arc.}
\label{fig:leaping}
\end{figure}

\end{proof}

Let~$n$ be the number of walls of the quasi wall system~$\W$ on the disk~$D$. From now on, changing the parameterization of the boundary circle $S^1=\partial D$, we assume that~$S^1$ is a circle of length~$n$, thus $S^1 = \R/n\Z$, and that the endpoints of the walls of~$\W$ are located at the semi integer points.
(This implies that the distance between two adjacent integer points is equal to~$1$.)
Therefore, on the universal cover of the cylinder $\C = D\setminus \{O\}$, which is the upper half plane~$\HH$, we have $\partial\W_\HH= \Z + \frac 12 \subseteq\R = \partial\HH$.

Note that the quasi wall system~$\W_\HH$ is periodic of period~$n$ (where~$n$ is the number of walls of~$\W$) in the sense that it is invariant by the horizontal translation of length~$n$. However, the following result implies that~$\W_\HH$ is also periodic with period~$2r$, where~$r$ is the number of strands of~$\W_\HH$; see Definition~\ref{def:strand}.

\begin{lemma} \label{lem:2r}
The sum of the widths of two adjacent arcs~$\alpha_0,\,\alpha_1$ of~${\W_\HH}$ is equal to~$2r$.
\end{lemma}

\begin{proof} Consider two adjacent arcs $\alpha_0=\widehat{a_0a_1}$ and $\alpha_1=\widehat{a_1a_2}$ as in Lemma~\ref{lem:zerotwo}. According to that lemma, each of the~$r$ strands of~$\W_\HH$ has exactly two stops in the interval $[a_0,a_2)$. Therefore there are~$2r$ semi-integers in that interval. It follows that $a_2-a_0=2r$.
\end{proof}

Denote by~$\strip_{[t,t+2r)} = [t,t+2r)\times[0,+\infty)$ a strip of width~$2r$ of the half-plane~$\HH$. The following result describes the arcs of the quasi wall system~$\W_\HH$ that are contained in such a strip.

\begin{proposition}\label{thm:strip_content}\mbox{}
\begin{enumerate}
\item\label{item:one_arc} Each strip $\strip_{[t,t+2r)}$ contains exactly one arc of each strand (and each of these arcs determines its strand completely).
\item\label{item:no_inter} The~$r$ arcs contained in a strip $\strip_{[t,t+2r)}$ do not intersect each other.
\item\label{item:cross_twice} Any pair of strands intersects each other exactly twice in the strip $\strip_{[t,t+2r)}$.
\end{enumerate}
\end{proposition}

\begin{proof}\mbox{}
\\\noindent\eqref{item:one_arc}
Consider a strand~$(\alpha_i)_{i\in\Z}$, with $\alpha_i=\widehat{a_ia_{i+1}}$. According to Lemma~\ref{lem:2r}, we have the equation $a_{i+2}=a_i+2r$ for all~$i$. This implies that the strip $\strip_{[t,t+2r)}$ contains exactly two stops and thus exactly one arc of the strand~$(\alpha_i)_i$. The same equation implies that two consecutive stops determine the strand.
\\\noindent\eqref{item:no_inter}
Consider a second strand~$(\beta_j)_{j\in\Z}$, with~$\beta_j=\widehat{b_jb_{j+1}}$. Assuming that two arcs~$\alpha_0$ and~$\beta_0$ of~$\W_\HH$ intersect, we want to show that they are not contained in a strip $\strip_{[t,t+2r)}$. We may assume without loss of generality that~$a_0<b_0$, therefore~$b_0\in (a_0,a_1)$. Since the strand~$(\beta_j)_j$ has a stop in the interval~$[a_0,a_1)$ by Lemma~\ref{lem:zerotwo} it cannot have a stop in~$[a_1,a_2)$. It follows that~$b_1>a_2=a_0+2r$, hence the arcs~$\alpha_0 = \widehat{a_0,a_1}$ and~$\beta_0=\widehat{b_0,b_1}$ are not contained in a strip of width~$2r$.
\\\noindent\eqref{item:cross_twice}
Consider two strands~$(\alpha_i)_{i\in\Z}$ and~$(\beta_j)_{j\in\Z}$ as above. 
Since $a_{i+2}=a_i+2$ as shown in~\eqref{item:one_arc}, the stand~$(\alpha_i)_i$ is invariant by the horizontal translation of displacement~$2r$.
The same holds with~$(\beta_j)_j$.
We want to show that they cross exactly twice in a strip~$\strip_{[t,t+2r)}$. 
By invariance under the horizontal translation of length~$2r$, we may choose~$t$ arbitrarily. 
For instance, we can choose~$t=a_0$. 
By Lemma~\ref{lem:zerotwo}, the strand~$(\beta_j)_j$ has stops in exactly one of the intervals~$(a_0,a_1)$ and~$(a_1,a_2)$. Thus, it intersects (twice) exactly one of the arcs~$\alpha_0=\widehat{a_0a_1}$,~$\alpha_1=\widehat{a_1a_2}$.
\end{proof}

We also note the following.

\begin{lemma}\label{lem:width_one}
In the quasi wall system~$\W_\HH$, there is an arc of width 1.
\end{lemma}

\begin{proof} Let $\alpha_0=\widehat{a_0a_1}$ be an arc that is minimal with respect to covering (\ie, $\alpha_0$ does not cover any arc of~$\W_\HH$). We want to show that $a_1-a_0=1$. By Lemma~\ref{lem:zerotwo}, each strand other than the one generated by~$\alpha_0$ has two stops in the interval $(a_0,a_2)$, both contained either in $(a_0,a_1)$ or in $(a_1,a_2)$. Thus, if the interval $(a_0,a_1)$ has any stop, it has in fact two stops of a strand, and therefore there is an arc of~$\W_\HH$ covered by~$\alpha_0$. However, this possibility is excluded by the minimality of~$\alpha_0$. Therefore, the interval $(a_0,a_1)$ has no stops and hence its endpoints~$a_0$ and~$a_1$ are consecutive semi-integers.
\end{proof}

\section{Proof of the discrete area lower bound} \label{sec:proof}

We can now proceed to the proof of the discrete area lower bound for simple discrete metric disks, see Theorem~\ref{theo:discrete-area}, making use of the previous notations and constructions.
Namely, let us prove the following.

\begin{theorem} \label{theo:discrete-area2}
The discrete area of every simple discrete metric disk of radius~$r$ is at least~$\frac 32 r^2$.
\end{theorem}

\begin{proof} Let~$(D,\W')$ be a simple discrete metric disk of radius~$r$ and center~$O$ that has minimal area. Recall that the punctured disk~$D\setminus\{O\}$ is identified with the flat cylinder $\C = S^1\times [0,+\infty)$. As shown in Section~\ref{sec:semi-circle}, $\W'$ is homotopic in~$\C$ to a quasi wall system~$\W$ made of standard arcs, such that~$(D,\W)$ is also a discrete disk of radius~$r$ centered at~$O$ and has the same area as~$(D,\W')$. Thus we must show that $\area(D,\W) \geq \frac 32 r^2$. Also, we may assume that the lift of $\W$ to the universal cover $\HH=\R\times[0,+\infty)$ is a quasi wall system~$\W_\HH$ such that $\partial\W_\HH = \Z + \frac 12 \subseteq \R = \partial\HH$ as in Section~\ref{sec:adjacent}.

Let~$t\in\R$ be a generic number. By Proposition~\ref{thm:strip_content}, the weighted number of self-intersections of the quasi wall system~$\W_\HH$ that lie in the strip $\strip_{[t,t+2r)}$ is
\begin{equation} \label{eq:area_strip}
\# I_{\Int \HH +\frac 12 \partial \HH}(W_\HH|_{S_{[t,t+2r)}})
= 2 \frac{r(r-1)}2 + \frac 12 2r
= r^2.
\end{equation}
The first term counts the crossings between the different strands: each pair of strands crosses twice, and the crossings are located in the interior of the half-plane~$\HH$. The second term counts, with weight~$\frac 12$, the intersections that lie in the boundary~$\partial\HH$; these are the intersections between adjacent arcs, that belong to the same strand. Thus, the discrete area of the disk~$(D,\W)$ is 
\[ \area(D,\W)
= \# I_{\Int \HH +\frac 12 \partial \HH}(W_\HH|_{S_{[t,t+n)}})
= \frac n{2r} \# I_{\Int \HH+\frac 12 \partial \HH}(W_\HH|_{S_{[t,t+2r)}})
= \frac{n}{2r}\,r^2, \]
where~$n$ is the number of walls of~$\W$. 

To finish we will show that $n\geq 3r$. Let $(\alpha_i=\widehat{a_ia_{i+1}})_{i\in\Z}$ be a strand of~$\W_\HH$ such that $a_0-a_{-1}=1$. Such a strand exists by Lemma~\ref{lem:width_one}. Moreover, we may assume that~$a_0=\frac 12$ and~$a_{-1}=-\frac 12$.
The interval $(a_0,a_1)$ has width~$2r-1$ (by Lemma~\ref{lem:2r}) and contains~$2r-2$ semi-integers.

Each of these semi-integers is either the startpoint or the endpoint of one of the~$r-1$ arcs that are covered by~$\alpha_0$; see Proposition~\ref{thm:strip_content}. Let~$b_0$ be the rightmost of the~$r-1$ startpoints. Note that
\begin{equation}\label{eq:xckeut} b_0\geq a_0+(r-1). \end{equation}
This point~$b_0$ is a stop of a strand $(\beta_j=\widehat{b_jb_{j+1}})_{j\in\Z}$. The arc~$\beta_0$ is covered by~$\alpha_0$ and the arc $\beta_1=\widehat{b_1b_2}$ intersects the arc~$\alpha_0$. The arcs~$\alpha_0$ and~$\beta_1$ cannot extend over a whole fundamental domain $S_{[t,t+n)}$ of the universal cover, by the property~\eqref{item:no_two_cover} of Proposition~\ref{thm:simple_interval_system}. Therefore, $n > b_2 - a_0 $. On the other hand, by Lemma~\ref{lem:2r} and the inequality~\eqref{eq:xckeut}, we have
\[ b_2 = b_0+2r \geq a_0 + 3r - 1. \]
We conclude that $n > 3r - 1$, or, equivalently, $n\geq 3r$, as we had to prove.
\end{proof}

\section{Simple discrete metric disks of minimal area} \label{sec:extremal1}

In this section, we analyze the equality case of Theorem~\ref{theo:discrete-area2}.

\begin{proposition} \label{prop:isotopy}
For every positive integer~$r$, there is a simple discrete metric disk of radius~$r$ and area~$\frac{3}{2} r^2$.
It is unique up to isotopy of the disk with the center fixed.
\end{proposition}

\begin{figure}[htbp]
\vspace{-2mm}
\def\svgwidth{128mm}
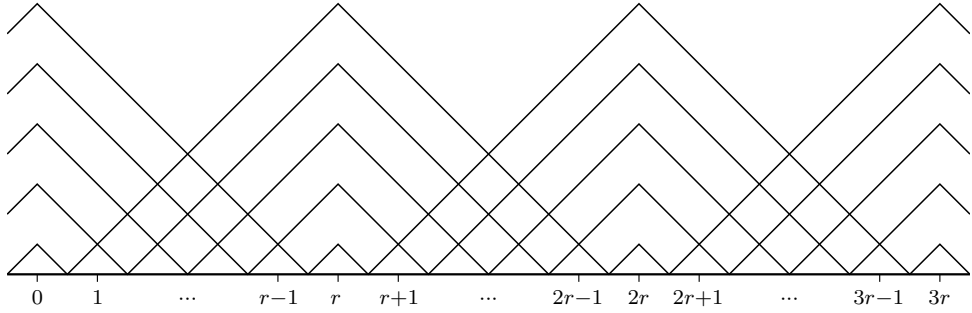
\caption{The lift~$W_\HH$ corresponding to an area minimizing simple discrete disk $(D,\W)$ of radius~$r=5$, where the quasi wall system~$\W$ consists of standard arcs.}
\label{fig:min_standard}
\end{figure}

\begin{figure}[htbp]
\def\svgwidth{50mm}
\begingroup%
  \makeatletter%
  \providecommand\color[2][]{%
    \errmessage{(Inkscape) Color is used for the text in Inkscape, but the package 'color.sty' is not loaded}%
    \renewcommand\color[2][]{}%
  }%
  \providecommand\transparent[1]{%
    \errmessage{(Inkscape) Transparency is used (non-zero) for the text in Inkscape, but the package 'transparent.sty' is not loaded}%
    \renewcommand\transparent[1]{}%
  }%
  \providecommand\rotatebox[2]{#2}%
  \newcommand*\fsize{\dimexpr\f@size pt\relax}%
  \newcommand*\lineheight[1]{\fontsize{\fsize}{#1\fsize}\selectfont}%
  \ifx\svgwidth\undefined%
    \setlength{\unitlength}{255.11811024bp}%
    \ifx\svgscale\undefined%
      \relax%
    \else%
      \setlength{\unitlength}{\unitlength * \real{\svgscale}}%
    \fi%
  \else%
    \setlength{\unitlength}{\svgwidth}%
  \fi%
  \global\let\svgwidth\undefined%
  \global\let\svgscale\undefined%
  \makeatother%
  \begin{picture}(1,1)%
    \lineheight{1}%
    \setlength\tabcolsep{0pt}%
    \put(0,0){\includegraphics[width=\unitlength,page=1]{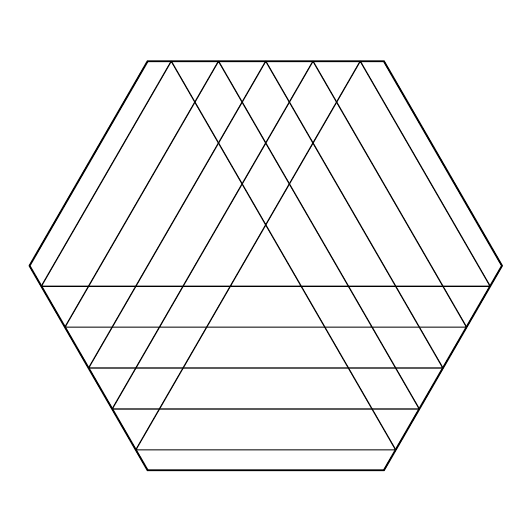}}%
    \put(0.48343391,0.48333702){\makebox(0,0)[lt]{\lineheight{1.25}\smash{\begin{tabular}[t]{l}$O$\end{tabular}}}}%
  \end{picture}%
\endgroup%

\caption{An area minimizing simple discrete disk $(D,\W)$ of radius $r=5$ where the topological disk~$D$ is an hexagon and the quasi wall system~$\W$ consists of straight lines.}
\label{fig:min_standard_hexagon}
\end{figure}

\begin{proof}
Recall the proof of Theorem~\ref{theo:discrete-area2}. 
Let~$\W'$ be a simple quasi wall system such that $(D,\W')$ is a simple discrete metric disk of radius~$r$ with minimal discrete area.
Consider the simple quasi wall system~$\W$ homotopic to~$\W'$ made of standard arcs.
To attain the lower bound on~$\area(D,\W)$ and so on~$\area(D,\W')$, we must have $n=3r$, therefore the inequality~\eqref{eq:xckeut} must be an equality. 
This implies that, for the~$r-1$ arcs covered by~$\alpha_0$, the~$r-1$ startpoints must precede the~$r-1$ endpoints in the interval $(a_0,a_1)$. 
In consequence, these~$r-1$ arcs together with the arc~$\alpha_0$ form a chain with respect to the covering relation; see Figure~\ref{fig:min_standard}. 
This implies that the~$r$ arcs are completely determined, and by Proposition~\ref{thm:strip_content}, so are the quasi wall systems~$\W_\HH$ and~$\W$, which are made of standard arcs. 
Thus, the quasi wall system~$\W_\HH$ contains all arcs of the form $\widehat{kr-s,kr+s}$ with~$k$ integer and $s\in(0,r)$ semi integer; see Figure~\ref{fig:min_standard}.
Similarly, the quasi wall system~$\W$ is obtained from~$\W_\HH$ by taking the quotient of~$\HH$ under the horizontal translation of length~$3r$; see Figure~\ref{fig:min_standard_hexagon}.
This proves the uniqueness of the simple discrete metric disk of minimal area, but only up to homotopy of the quasi wall system.
The uniqueness up to isotopy of the disk follows from the next result.

\begin{lemma}\label{thm:isotopy}
Let~$\W$ and~$\W'$ be two simple quasi wall systems on the disk~$D$, homotopic in~$D \setminus \{0\}$ and forming no triangle in~$D \setminus \{O\}$.
Then there is an isotopy of~$D$ which fixes~$O$ and carries~$\W$ to~$\W'$.
\end{lemma}

\begin{proof}
We proceed by induction in the number~$n$ of walls. 
The case $n=1$ is trivial.
In general, we argue as follows.

Let~$\gamma$ be a wall of~$\W$ that covers no other wall of~$\W$; see Definition~\ref{def:cover}.
The curve~$\gamma$ divides the disk~$D$ into two topological closed disks~$A$ and~$B$ which intersect along~$\gamma$, with~$O \in A$.
The part of~$\W$ that lies in~$B$ consists of $k\geq 0$ arcs going from~$\gamma$ to~$\partial B\setminus\gamma$. 
These arcs are pairwise disjoint, otherwise they would form a triangle in~$B\subseteq D\setminus\{O\}$. 
The part of~$\W$ that lies in~$A$, excluding~$\gamma$, is a quasi wall system on~$A$ with~$n-1$ walls.

Let~$\gamma'$ be the wall of~$\W'$ homotopic to~$\gamma$ in~$D\setminus\{O\}$.
We apply to~$\W'$ a first isotopy of~$D\setminus\{O\}$ to ensure that $\gamma'=\gamma$. 
The wall~$\gamma'$ does not cover any other wall~$\beta'$ of~$\W'$, otherwise the wall~$\beta$ of~$\W$ homotopic to~$\beta'$ would cross~$\gamma$ twice. 
Similarly as in~$\W$, the part of~$\W'$ lying in~$B$ consists of~$k$ pairwise disjoint arcs going from~$\gamma$ to~$\partial B\setminus\gamma$. 
Thus, by applying a second isotopy, we may ensure that $\W'\cap B=\W\cap B$. Finally, we get $(\W'\setminus\gamma')\cap A = (\W\setminus\gamma)\cap A$ by applying an isotopy of the disk~$A$ fixing~$O$, whose existence is guaranteed by the inductive hypothesis.
\end{proof}

Now, the walls of~$\W$ do not delimit any triangle in~$D \setminus \{O\}$ (where each side lies in a wall); see Figures~\ref{fig:min_standard_hexagon} and~\ref{fig:min_standard}.
Since two arcs of~$\W$ intersect each other if and only if the same holds with the corresponding homotopic arcs of~$\W'$, we deduce that the walls of~$\W'$ do not form any triangle in~$D \setminus \{O\}$ either.
The uniqueness of the simple discrete metric disk of minimal area up to isotopy of the disk fixing~$O$ follows from Lemma~\ref{thm:isotopy}.
\end{proof}

\begin{remark}
The isotopy between~$\W$ and~$\W'$ can also be derived from \cite{graaf1997making}, where it is proved that two wall systems on a closed surface which are homotopic to each other and are both in minimally crossing position (\ie, they attain the minimum number of self-intersections possible in their homotopy class) can be obtained one from the other by isotopies and triangle flip moves (called ``type III moves'' in~\cite{graaf1997making}). 
Strictly speaking, we first need to adapt this result to quasi wall systems on surfaces with boundary.
Since~$\W$ and~$\W'$ do not form any triangle in $D\setminus\{0\}$, we conclude that they are isotopic in~$D$.
\end{remark}

\section{Construction of almost minimizing Finsler disks}
\label{sec:extremal2}

In this section, we construct a Finsler disk of radius~$r$ with minimizing interior geodesics whose area is arbitrarily close to the lower bound~$\frac{6}{\pi} r^2$ given by Theorem~\ref{theo:6pi-intro}.

\medskip

Let us first go over Busemann's construction of projective metrics in relation with Hilbert's fourth problem.
We refer to \cite{busemann1976problem}, \cite{pogorelov1979hilbert}, \cite{alexander1978planes}, \cite{szabo1986hilbert}, \cite{papadopoulos2014hilbert} and references therein for an account on the subject.

The space~$\Gamma$ of oriented lines in~$\R^2$ can be identified with~$S^1 \times \R$.
Under this identification, an oriented line~$\gamma$ is represented by a pair~$(\e^{\ii\theta},p)$ where~$\e^{\ii\theta}$ is the direction of the oriented line~$\gamma$ and $p = \langle \overrightarrow{OH} \times \e^{\ii\theta},\overrightarrow{e_z} \rangle$ is the signed distance from the origin~$O$ to~$\gamma$.
Here,~$H$ is a point of~$\gamma$, the vector~$\overrightarrow{e_z}$ is the third vector in the canonical basis of~$\R^3$, thus it is a unit vector orthogonal to~$\R^2$, and ``$\times$'' is the vector product in~$\R^3$.

\begin{definition} \label{def:borel}
Let~$\mu$ be a (nonnegative) Borel measure on~$\Gamma$.
Consider the following conditions:
\begin{enumerate}
\item the measure is invariant under the involution of~$\Gamma$ reversing the orientation of lines; \label{adm1}
\item the measure of every compact subset of~$\Gamma$ is finite; \label{adm2}
\item the set of all oriented lines passing through any given point of~$\R^2$ has measure zero; \label{adm3}
\item the set of all oriented lines passing through any given line segment in~$\R^2$ has positive measure. \label{adm4}
\end{enumerate}
A Borel measure~$\mu$ satisfying \eqref{adm1}-\eqref{adm3} induces a length function
\[\length_\mu(\alpha)
= \frac{1}{4} \int_{\gamma \in \Gamma} \#(\gamma \cap \alpha) \, \diff\mu(\gamma)
\] defined for any curve~$\alpha$ in the plane~$\R^2$. For this kind of length function, straight segments are shortest paths, therefore the pseudo-distance associated to this length function is
\[ d_\mu(x,y)
= \frac 14 \, \int_{\gamma \in \Gamma} \#(\gamma \cap [x,y]) \, \diff\mu(\gamma)
= \frac 14 \mu(\Gamma_{[x,y]}), \]
where~$\Gamma_A$ denotes the set of lines $\gamma\in\Gamma$ that intersect a subset or point~$A$ contained in the plane~$\R^2$.
The pseudo-distance $d_\mu$ is \term{projective}, which means that $d(x,z)=d(x,y)+d(y,z)$ for every $x,y,z \in \R^2$ with $y \in [x,z]$, and in fact every continuous projective distance is obtained from a unique measure~$\mu$; see~\cite{alexander1978planes}.
If~$\mu$ also satisfies~\eqref{adm4} then~$d_\mu$ is a projective distance (and vice-versa).
\end{definition}

For example, the product measure~$\lambda$, given by $\diff\lambda = \diff\theta \diff p$, yields the Euclidean distance $d_\lambda(x,y)=|y-x|$.

\medskip

The projective distance induced by a Borel measure satisfying the conditions \eqref{adm1}--\eqref{adm4} is not Finsler in general.
Borel measures inducing a Finsler metric can be characterized as follows; see \cite{paiva2005symplectic} for a presentation of this result due to Pogorelov~\cite{pogorelov1979hilbert} and~\cite{paiva2010finsler} for a generalization.

\begin{theorem} \label{thm:projective_finsler}
Let~$\mu$ be a Borel measure on~$\Gamma$ satisfying \eqref{adm1}--\eqref{adm4}.
The distance~$d_\mu$ is Finsler if and only if~$\mu$ is a positive smooth measure.
In this case, the smooth measure on~$\Gamma$ induced by the symplectic form associated to the Finsler metric, see~\eqref{eq:muG}, coincides with~$\mu$.
\end{theorem}

Here, a measure~$\mu$ on~$\Gamma$ is \term{(positive) smooth} if it admits a (positive) smooth function~$h$ as density, that is, $\diff\mu = h \diff\lambda$.

\begin{remark}\label{rmk:projective_geodesics} The geodesics of a plane with a projective Finsler metric~$d_\mu$ are the straight lines parametrized by $\mu$-length.
Therefore, a plane with a projective Finsler metric has minimizing geodesics.
\end{remark}

We may define the area of a Borel set~$D$ in the plane with a measure~$\mu$ on~$\Gamma$ satisfying \eqref{adm1}--\eqref{adm3} by the Santal\'o+Blaschke formula~\eqref{eq:croftonsantalo}
\begin{equation}
\label{eq:santaloblaschke_mu}
\area_\mu(D) 
= \frac 1{8\pi} \int_{\gamma_0\in\Gamma} \int_{\gamma_1\in\Gamma}
\#(\gamma_0\cap\gamma_1\cap D) \, \diff\mu(\gamma_1) \diff\mu(\gamma_0).
\end{equation}
In other terms, the area measure is the normalized pushforward measure
\begin{equation}
\label{eq:area_mu}
\area_\mu = \frac 1{8\pi} i_*(\mu\times\mu)
\end{equation}
where $i:\Gamma\times\Gamma\setminus \Delta_\Gamma \to \R P^2$ maps each ordered pair of different lines to its intersection point in the projective plane $\R P^2\supseteq \R ^2$. (Note that the diagonal $\Delta_\Gamma$ has measure zero because~$\mu$ has no atoms.)
This area function coincides with Holmes--Thompson area if the metric is Finsler; see~\eqref{eq:croftonsantalo}.

\subsection{Construction of a non-Finsler extremal disk}
\begin{figure}[htbp] 
\vspace{0.8cm}
\def\svgwidth{6cm}
\begingroup%
  \makeatletter%
  \providecommand\color[2][]{%
    \errmessage{(Inkscape) Color is used for the text in Inkscape, but the package 'color.sty' is not loaded}%
    \renewcommand\color[2][]{}%
  }%
  \providecommand\transparent[1]{%
    \errmessage{(Inkscape) Transparency is used (non-zero) for the text in Inkscape, but the package 'transparent.sty' is not loaded}%
    \renewcommand\transparent[1]{}%
  }%
  \providecommand\rotatebox[2]{#2}%
  \newcommand*\fsize{\dimexpr\f@size pt\relax}%
  \newcommand*\lineheight[1]{\fontsize{\fsize}{#1\fsize}\selectfont}%
  \ifx\svgwidth\undefined%
    \setlength{\unitlength}{425.19685039bp}%
    \ifx\svgscale\undefined%
      \relax%
    \else%
      \setlength{\unitlength}{\unitlength * \real{\svgscale}}%
    \fi%
  \else%
    \setlength{\unitlength}{\svgwidth}%
  \fi%
  \global\let\svgwidth\undefined%
  \global\let\svgscale\undefined%
  \makeatother%
  \begin{picture}(1,0.63333333)%
    \lineheight{1}%
    \setlength\tabcolsep{0pt}%
    \put(0.49451886,0.74198907){\color[rgb]{0,0,0}\makebox(0,0)[lt]{\lineheight{0}\smash{\begin{tabular}[t]{l} \end{tabular}}}}%
    \put(0.40569444,0.22479167){\color[rgb]{0,0,0}\makebox(0,0)[lt]{\lineheight{0}\smash{\begin{tabular}[t]{l} \end{tabular}}}}%
    \put(-0.21166667,0.7275){\color[rgb]{0,0,0}\makebox(0,0)[lt]{\lineheight{0}\smash{\begin{tabular}[t]{l} \end{tabular}}}}%
    \put(0,0){\includegraphics[width=\unitlength,page=1]{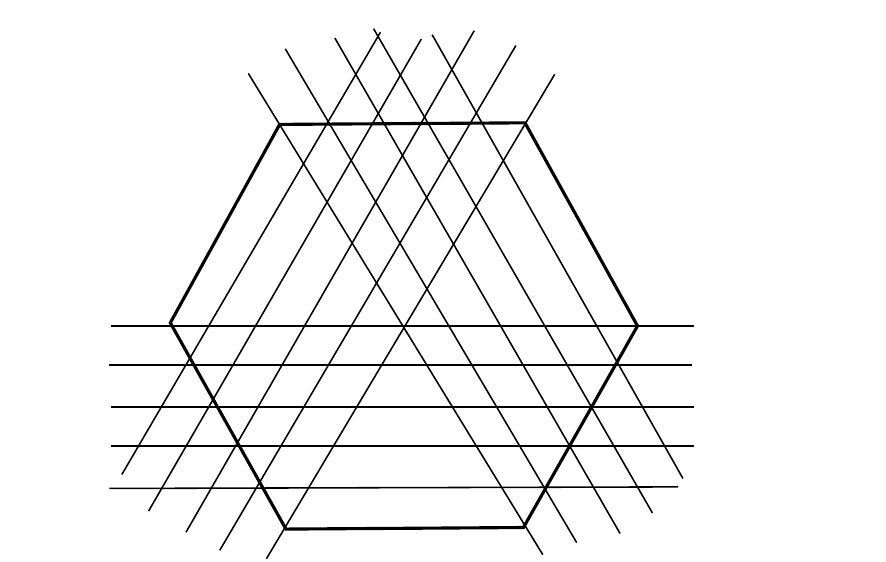}}%
    \put(0.00352778,0.17011111){\color[rgb]{0,0,0}\makebox(0,0)[lt]{\lineheight{0}\smash{\begin{tabular}[t]{l}$L_0$\end{tabular}}}}%
    \put(0.22048611,0.62519445){\color[rgb]{0,0,0}\makebox(0,0)[lt]{\lineheight{0}\smash{\begin{tabular}[t]{l}$L_1$\end{tabular}}}}%
    \put(0.61206944,0.62519445){\color[rgb]{0,0,0}\makebox(0,0)[lt]{\lineheight{0}\smash{\begin{tabular}[t]{l}$L_2$\end{tabular}}}}%
  \end{picture}%
\endgroup%
 
\caption{Extremal pseudo-metric disk} \label{fig:hexagon}
\end{figure}
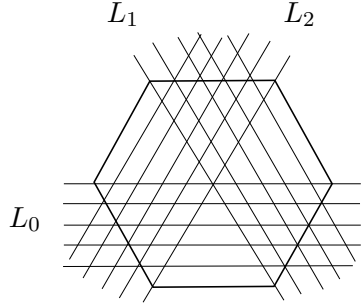

Let us construct a non-Finsler projective pseudo-metric disk satisfying the equality case in Theorem~\ref{theo:6pi-intro}.
Consider the three pairs of one-parameter families~$L_k^\pm$ of oriented lines in~$\R^2$ defined as
\[
\begin{array}{cccc}
L_k^+: & \R_+ & \to & \Gamma = S^1 \times \R \\
 & t & \mapsto & (\e^{\ii \frac{2k\pi}{3}},t)
\end{array}
\quad \mbox{ and } \quad
\begin{array}{cccc}
L_k^-: & \R_+ & \to & \Gamma = S^1 \times \R \\
 & t & \mapsto & (\e^{\ii (\frac{2k\pi}{3}+\pi)},-t)
\end{array}
\]
where $k \in \{0,1,2\}$; see Figure~\ref{fig:hexagon}.
Note that the lines~$L_k^+(t)$ and~$L_k^-(t)$ only differ by their orientation.
We will sometimes denote these families of lines by~$L_k$ when the orientation does not matter.
Consider the (nonsmooth) Borel measure on~$\Gamma$
\[ \mu_{\ext} = \nu_0+\nu_1+\nu_2 \]
where 
\[ \nu_k = \tfrac{1}{2} [(L_k^+)_*(\mathcal{L}) + (L_k^-)_*(\mathcal{L})]  \]
is the average of the push-forwards to~$\Gamma$ of the Lebesgue measure~$\mathcal{L}$ on~$\R_+$.
Let~$D_k$ be the line passing through~$O$ orthogonal to~$L_k$.
Let $\overline{D}_k \subseteq D_k$ be the ray from~$O$ that intersects orthogonally every line~$L_k(t)$.
Denote by~$\pi_k$ the orthogonal projection of~$\R^2$ to~$D_k$.
By construction, the $d_{\nu_k}$-pseudo-distance between two points $x,y \in \R^2$ is equal to one quarter times the Euclidean length of the projection of~$[x,y]$ to~$D_k$ lying in~$\overline{D}_k$. That is,
\[ d_{\nu_k}(x,y)
= \frac 14 \length(\pi_k([x,y]) \cap \overline{D}_k) \leq \frac 14|x-y| \]
for every $x,y \in \R^2$.
Furthermore,
\[ d_{\mu_{\ext}}(x,y)
= \sum_{k=0,1,2} d_{\nu_k}(x,y)
= \frac 14 \sum_{k=0,1,2} \length(\pi_k([x,y]) \cap \overline{D}_k). \]
Observe also that the measure~$\mu_{\ext}$ satisfies \eqref{adm1}-\eqref{adm3}, but not~\eqref{adm4}.
Thus, $d_{\mu_{\ext}}$ is a projective pseudo-distance on~$\R^2$.

The disk~$D_{\mu_{\ext}}(r)$ of radius~$r$ for the pseudo-distance~$d_{\mu_{\ext}}$ with center the origin~$O$ of~$\R^2$ is the minimal regular hexagon containing the Euclidean disk of radius~$4r$, whose vertices are $\frac{2\sqrt{3}}{3} 4r \e^{\ii \, k \frac{\pi}{3}}$ for $k \in \{0,\dots,5\}$; see Figure~\ref{fig:hexagon}. 
A direct computation using~\eqref{eq:santaloblaschke_mu} shows that its area is~$\frac 6\pi r^2$.
Thus, the disk~$D_{\mu_{\ext}}(r)$ is a non-Finsler projective pseudo-metric disk satisfying the equality case in Theorem~\ref{theo:6pi-intro}. 
One can think of it as an extremal (degenerate) metric for the problem considered.
Observe also that~$D_{\mu_{\ext}}(r)$ is not rotationally symmetric.

\begin{remark} \label{rem:extremal}
By identifiying all pairs of points at zero pseudo-distance, the pseudo-metric disk~$D_{\mu_{\ext}}(r)$ identifies with the closed ball~$D(r)$ of radius~$r$ centered at the tip of a cone composed of three copies of a quadrant of the $\ell^1$-plane glued together.
It follows from a direct computation that the Holmes--Thompson area of the disk~$D(r)$ is equal to~$\frac{6}{\pi}r^2$.
Defined in this way, the metric on~$D(r)$ is non-Finsler (\eg, it has a singularity at the origin and the tangent norms are neither smooth nor strongly convex) but can still be thought of as an extremal (degenerate) metric.
Note that the (pseudo)-metrics on~$D_{\mu_{\ext}}(r)$ and~$D(r)$ can be viewed as continuous versions of the extremal simple discrete disk; see~Section~\ref{sec:extremal1}.
\end{remark}


\subsection{Construction of a Finsler nearly extremal disk}
In the rest of this section, we explain how to modify the pseudo-metric~$d_{\mu_{\ext}}$ so as to obtain a projective Finsler disk of radius~$r$ whose area is arbitrarily close to~$\frac{6}{\pi} r^2$.
First, the projective pseudo-metric~$d_{\mu_{\ext}}$ can be approximated by a projective metric by simply adding to~$\mu$ a multiple~$\eps\lambda$ of the uniform measure~$\lambda$ (given by $\diff\lambda = \diff\theta \diff p$) so that the point~\eqref{adm4} is also satisfied; this changes~$d_{\mu_{\ext}}$ by adding~$\eps$ times the Euclidean distance.
This projective metric is not Finsler, but in turn it can be approximated by a Finsler metric; see~\cite{pogorelov1979hilbert}.
More generally, every projective distance~$d_\mu$, where~$\mu$ is a Borel measure satisfying \eqref{adm1}-\eqref{adm4}, can be approximated by a projective Finsler distance on every compact set of~$\R^2$.
Thus, by Theorem~\ref{thm:projective_finsler}, there exists a sequence~$\mu_n$ of positive smooth measures on~$\Gamma$ such that the corresponding sequence of Finsler distances~$d_{\mu_n}$ uniformly converges to~$d_{\mu_{\ext}}$ on every compact set of~$\R^2$.
This approximation result is obtained by a convolution argument on the distance function~$d_\mu$.
Although it is possible that the measures~$\mu_n$ weakly converge to~$\mu_{\ext}$, this issue is not addressed in~\cite{pogorelov1979hilbert}.
This leads us to consider a slightly different approach.
Instead of regularizing the distance function, we smooth out the measure~$\mu_{\ext}$ and show that the corresponding projective Finsler distance converges to~$d_{\mu_{\ext}}$.
This alternative approach to the regularization of a projective distance provides a weak convergence of measure by construction, which allows us to estimate areas as well as distances.

\medskip

We proceed as follows. First, we truncate the measure~$\mu_{\ext}$ by setting a bound for the absolute value of the~$p$ coordinate of the lines $\gamma\in\Gamma$. In this way, we obtain a probability measure~$\mu_0$ on~$\Gamma$, without changing the corresponding distance function in a neighborhood of the origin. Similarly, we truncate the uniform measure~$\lambda$ to a probability measure~$\lambda_0$. This enables us to use standard theorems on weak convergence of probability measures.

Let us now describe the convolution process.
For $\eps>0$, let~$h_\eps$ be a smooth nonnegative function on $\Gamma = \R/2\pi\Z \times \R$, with support in~$(-\eps,\eps) \times (-\eps,\eps)$, such that $\int_\Gamma h_\eps(\theta,p) \, \diff\theta \diff p = 1$.
For each~$\eps>0$, consider the positive smooth measure~$\mu_\eps$ on~$\Gamma$ with density $h_\eps \convol \mu_0$, that is,
\[ \diff\mu_\eps = (h_\eps \convol \mu_0) \, \diff\lambda \]
where $h_\eps \convol \mu_0$ is the smooth function on~$\Gamma$ defined by the convolution
\[ h_\eps \convol \mu_0(\gamma)
= \int_\Gamma h_\eps(\gamma - \gamma') \, \diff\mu_0(\gamma') \]
and~$\lambda$ is the standard product measure on $\Gamma = \R/2\pi\Z \times \R$, given by $\diff\lambda = \diff\theta \diff p$.
By~\cite[\S1.4.3]{bogachev2018weak}, the smooth measure~$\mu_\eps$ weakly converges to~$\mu_0$ as~$\eps$ goes to zero.
Define also the measure
\[ \mu_\eps^+ = (1-\eps) \mu_\eps + \eps \, \lambda_0, \]
which also converges to~$\mu_0$ as~$\eps\to 0$.
By Theorem~\ref{thm:projective_finsler}, the distance~$d_{\mu_\eps^+}$ induced by~$\mu_\eps^+$ is a projective Finsler distance on a neighborhood of the origin in~$\R^2$.

\medskip

To approximate distances and areas, we have the following tools.

\begin{lemma} \label{thm:unifapprox}
Let~$\mu$ and~$\mu_n$ be probability measures on~$\Gamma$ satisfying the conditions \eqref{adm1}--\eqref{adm3} of Definition~\ref{def:borel}.
If~$\mu_n$ weakly converges to~$\mu$, then the distance~$d_{\mu_n}$ converges uniformly to~$d_\mu$ on every compact subset of~$\R^2$.
\end{lemma}

\begin{proof}
Note first that the distance between two points $x,y\in \R^2$ is
\[ d_\mu(x,y) = \mu(\Gamma_{[x,y]}) \]
where~$\Gamma_{[x,y]}$ denotes the set of lines that intersect the segment~$[x,y]$.
Thus, for a specific pair of points $x,y$, the weak convergence $\mu_n \to \mu$ implies that $d_{\mu_n}(x,y) \to d_\mu(x,y)$ by the portmanteau theorem~\cite[Theorem 2.1]{billingsley1968convergence}, since $\Gamma_{[x,y]}$ is a continuity set for $\mu$.
That is, its boundary
\[ \partial \Gamma_{[x,y]} = \Gamma_x \cup \Gamma_y \] 
(where ~$\Gamma_z$ is the set of lines that contain a point~$z$)
has measure $\mu(\partial \Gamma_{[x,y]}) = 0$ since $\mu(\Gamma_z) = 0$ for each point $z \in \R^2$ by condition~\eqref{adm3} on~$\mu$.

To show that this convergence holds uniformly for $x$, $y$ in any given compact set $K\subseteq\R^2$, let~$\mathcal A$ be the family of sets $\Gamma_{[x,y]}$ for ~$x, y \in K$.
According to Theorems~2 and~3 from \cite{billingsley1967uniformity}, to show uniform convergence $\mu_n(A)\to\mu(A)$ for all sets $A\in\mathcal A$, it is sufficient to show that $\mu(B_\delta(\partial A))\to 0$ uniformly as $\delta\to 0$, where~$B_\delta(S)$ denotes the $\delta$-neighborhood of a set $S\subseteq\Gamma$ (say, with respect to the supremum distance in terms of the coordinates $\theta,\,p$).
Now, 
\[B_\delta(\partial \Gamma_{[x,y]}) = B_\delta(\Gamma_x) \cup B_\delta(\Gamma_y),\]
therefore it suffices to show that $\mu(B_\delta(\Gamma_x))\to 0$ uniformly for all $x\in K$ as $\delta\to 0$.
Suppose that this is not the case.
Then there are sequences $\delta_m\to 0$ and $x_m\to x\in K$ such that $\mu(B_{\delta_m}(\Gamma_{x_m}))$ does not tend to zero.
However, we also have
\[ B_{\delta_m}(\Gamma_{x_m}) \subseteq B_{\delta_m'}(\Gamma_x) \]
for some sequence $\delta_m'\to 0$ (namely, $ \delta_m' = \delta_m + |x_m-x| $), which yields $B_{\delta_m'}(\Gamma_x) \to \Gamma_x$ and therefore
\[ \mu(B_{\delta_m}(\Gamma_{x_m})) \leq \mu(B_{\delta_m'}(\Gamma_x)) \to \mu(\Gamma_x) = 0. \]
This contradiction finishes the proof.
\end{proof}

\begin{lemma}\label{thm:area_convergence} 
Let~$\mu$ and~$\mu_n$ be probability measures on~$\Gamma$ satisfying the conditions \eqref{adm1}--\eqref{adm3} of Definition~\ref{def:borel}.
If~$\mu_n$ weakly converges to~$\mu$, then $\area_{\mu_n}(K)$ converges to $\area_\mu(K)$ for every compact set $K\subseteq\R^2$ such that $\mu(\partial K)=0$.
\end{lemma}

\begin{proof} We will use some properties of weak convergence of probability measures; see \cite{billingsley1968convergence}.
Since~$\mu_n$ weakly converges to~$\mu$, it follows from~\cite[Example 3.2]{billingsley1968convergence} that the product measure $\mu_n\times\mu_n$ converges weakly to $\mu\times\mu$ on $\Gamma\times\Gamma$.
Restricting to the set $\Gamma\times\Gamma\setminus\Delta_\Gamma$, the measures $\mu\times\mu$ and $\mu_n\times \mu_n$  are still probability measures since the diagonal $\Delta_\Gamma\subseteq\Gamma\times\Gamma$ has zero measure because $\mu$ and $\mu_n$ have no atoms.
Moreover, since the diagonal~$\Delta_\Gamma$ is a closed set%
, the product measure~$\mu_n\times\mu_n$ weakly converges to~$\mu\times\mu$ on~$\Gamma\times\Gamma\setminus\Delta_\Gamma$ by condition (iv) of the portmanteau theorem~\cite[Theorem 2.1]{billingsley1968convergence}. Furthermore, since the function $i:\Gamma\times\Gamma\setminus\Delta_\Gamma \to \R P^2$ is continuous, the pushforward measure $i_*(\mu_n\times\mu_n)$ weakly converges to $i_*(\mu\times\mu)$ by the definition of weak convergence; see~\cite[p. 14]{billingsley1968convergence}. Therefore, the area measure $\area_{\mu_n}=\frac 1{8\pi}i_*(\mu_n\times\mu_n)$ weakly converges to $\area_{\mu}$, with both area measures considered as probability measures on the projective plane~$\R P^2$; see~\eqref{eq:area_mu}. Finally, to show that $\area_{\mu_n}(K)\to\area_{\mu}(K)$, we must check, according to part (v) of the portmanteau theorem~\cite[Theorem 2.1]{billingsley1968convergence}, that~$K$ is a continuity set of $\area_{\mu}$, which by definition means that $\mu(\partial_{\R P^2}K)=0$. This follows from the facts that~$K$ is compact and $\mu(\partial_{\R^2}K)=0$.
\end{proof}

Consider the disk~$D_{\mu_\eps^+}(r)$ centered at~$O$ of radius~$r$ for the distance~$d_{\mu_\eps^+}$. The numbers $r>0$ and $\eps>0$ are small enough so that the truncations of~$\mu_{\ext}$ and~$\lambda$ have no effect on the disk~$D_{\mu_\eps^+}(r)$. The number~$r$ is fixed while~$\eps$ goes to 0.

\begin{proposition}\label{thm:nearly_extremal_disks}
The disk $D_{\mu_\eps^+}(r)$ is a projective Finsler disk with minimizing interior geodesics, whose area converges to~$\frac{6}{\pi} r^2$, as~$\eps$ goes to zero. Therefore, the area lower bound in Theorem~\ref{theo:6pi-intro} is sharp.
\end{proposition}

\begin{proof} The fact that~$d_{\mu_\eps^+}$ is a projective Finsler metric follows from Theorem~\ref{thm:projective_finsler} and the fact that its geodesics are minimizing was stated in Remark~\ref{rmk:projective_geodesics}.

To compute the area of the disk $D_{\mu_\eps^+}(r)$ we proceed as follows.
By uniform convergence of the metrics (Lemma~\ref{thm:unifapprox}), for every~$\delta > 0$ and every $\eps>0$ small enough, we have
\[ 
D_{\mu_0}(r-\delta) \subseteq D_{\mu_\eps^+}(r) \subseteq D_{\mu_0}(r+\delta).
\]
Therefore,
\[ \area_{\mu_\eps^+}(D_{\mu_0}(r-\delta))
  \leq \area_{\mu_\eps^+}(D_{\mu_\eps^+}(r))
  \leq \area_{\mu_\eps^+}(D_{\mu_0}(r+\delta)). \]
Since the sets $D_{\mu_0}(r\pm\delta)$ are compact and have boundary of $\mu_0$-measure zero, Lemma~\ref{thm:area_convergence} shows that
\[ \area_{\mu_\eps^+}(D_{\mu_0}(r\pm\delta)) \to
   \area_{\mu_0}(D_{\mu_0}(r\pm\delta)) \]
as $\eps\to 0$.
Since this holds for each $\delta>0$, we conclude that
\[ \area_{\mu_\eps^+}(D_{\mu_\eps^+}(r)) \to \area_{\mu_0}(D_{\mu_0}(r)) \]
as $\eps\to 0$.
\end{proof}

\section{Appendix: Differentiability of distance-realizing paths on Finsler surfaces with boundary} \label{sec:regularity}

Consider a smooth manifold~$M$ with smooth boundary endowed with a Finsler metric~$F$.
Recall that a distance-realizing curve is a curve~$\alpha:I \to M$ defined on an interval~$I \subseteq \R$ such that 
\[ d_F(\alpha(t),\alpha(t')) = t' - t \]
for every~$t<t'$.

If the manifold has an empty boundary (or, more generally, a convex boundary), then its distance-realizing curves satisfy a differential equation, and it is therefore clear that they are smooth.
However, if the boundary is not convex, then the distance-realizing curves are not $C^2$ in general, and they are not even determined by their initial velocity vector.
This happens, for instance, on the Euclidean plane minus an open disk.

In the case of Riemannian manifolds with boundary, it was claimed in \cite{wolter1979interior} and~\cite{alexander1981geodesics} that distance-realizing curves are~$C^1$. 
This result can also be recovered from~\cite{lytchak2006holder} by gluing together two copies of a Riemannian manifold~$M$ along their boundaries.
The Riemannian metric obtained on the resulting double manifold~$N$ is $\alpha$-H\"older continuous for any $\alpha \in (0,1]$; see~\cite[Example~3.3]{lytchak2006holder}.
By~\cite{lytchak2006holder}, the geodesics on~$N$ are~$C^1$ (and even~$C^{1,\frac{\alpha}{2-\alpha}}$), from which we can deduce that the distance-realizing curves on~$M$ are also~$C^1$.
This argument does not hold for Finsler metrics. 
Indeed, the double of a Finsler metric is not even a continuous Finsler metric in general.

\medskip

Here, by adapting the argument of~\cite{alexander1981geodesics}, we prove that the same result holds for Finsler surfaces.

\begin{theorem} \label{thm:shortest_are_C1}
On a Finsler surface~$M$ with boundary, every distance-realizing curve~$\alpha:I\to M$ is~$C^1$.
Furthermore, the velocity vectors~$\alpha'(t)$ have unit norm.
\end{theorem}

Let us introduce some technical definitions.
We assume without loss of generality that the surface~$M$ is the closed upper half of~$\R^2$.

\begin{definition}
Let~$\alpha:I\to M$ be a continuous curve, where $I\subseteq\R$ is an interval.
Fix~$t_0\in I$ and denote~$x_0=\alpha(t_0)$.
An \term{arrival velocity} of~$\alpha$ at~$t_0$ is a vector~$v\in T_{x_0}M$ that is an accumulation point of the set of vectors 
\[ V^-=
\left\{ \frac{\alpha(t)-\alpha(t_0)}{t-t_0} \mid t<t_0 \right\} \]
as~$t$ goes to~$t_0$.
Similarly, a \term{departure velocity} of~$\alpha$ at~$t_0$ is a vector~$v\in T_{x_0}M$ that is an accumulation point of the set of vectors
\[ V^+ = \left\{ \frac{\alpha(t)-\alpha(t_0)}{t-t_0} \mid t>t_0 \right\} \]
as~$t$ goes to~$t_0$.
Note that if~$\alpha$ is differentiable on the left (resp. right) at~$t_0$, then~$\alpha$ has exactly one arrival (resp. departure) velocity at~$t_0$.
\end{definition}

We begin by proving a weak differentiability result.

\begin{lemma} \label{lem:right}
Let $(M,F)$ be a Finsler manifold with boundary and let~$\alpha:I\to M$ be a distance-realizing curve.
Fix~$t_0\in I$ and denote~$x_0=\alpha(t_0)$.
Then
\begin{enumerate}
\item The curve~$\alpha$ has at least one arrival velocity and one departure velocity at~$t_0$ (unless~$t_0=\min I$ or~$t_0=\max I$, respectively).
\item Every arrival or departure velocity~$v$ has norm~$F_{x_0}(v)=1$.
\item If the curve~$\alpha$ is differentiable on one side at an interior point~$t_0$ of~$I$, then~$\alpha$ is differentiable at~$t_0$.
\end{enumerate}
\end{lemma}

\begin{proof} 
By continuity of the Finsler metric at~$x_0$, we can bound~$F_x$ below and above by two multiples of the norm~$F_{x_0} = | \cdot |$ for every~$x$ close enough to~$x_0$.
That is,
\[ \lambda^-\,|v| \leq F_x(v) \leq \lambda^+\,|v| \]
for every~$v\in\R^n$, which in turn implies that
\[ \lambda^-\,|x-x_0| \leq d_F(x_0,x) \leq \lambda^+\,|x-x_0|. \]
This implies that the sets of vectors~$V^\pm$ are bounded when~$t$ goes to~$t_0$, which implies the first claim. 
In fact, as~$x$ goes to~$x_0$, the optimal coefficients~$\lambda^\pm$ converge to~$1$, which implies the second claim.

To prove the last claim, we assume that the curve~$\alpha$ is differentiable on the left at an interior point~$t_0$ of~$I$.
(The argument is similar if~$\alpha$ is differentiable on the right at~$t_0$.)
Let~$v^-$ be the arrival tangent vector. 
Let us prove that~$\alpha$ is differentiable on the right at~$t_0$ and has departure tangent vector~$v^+=v^-$. 
By contradiction, assume that the set of vectors~$V^+$ has an accumulation point~$v^+\neq v^-$ as~$t$ goes to~$t_0$. 
As already noticed in the second claim, we have~$|v^-|=|v^+|=1$.
Since the norm~$F_{x_0}= | \cdot |$ is strictly convex, we also have~$|v^-+v^+|<2$.
Let~$\tau_m\to 0$ be a decreasing sequence of positive numbers such that
\[ y_m^+ = \alpha(t_0+\tau_m) = \alpha(t_0) + \tau_m v^+ + o(\tau_m). \]
Since~$\alpha$ is differentiable on the left at~$t_0$, we also have
\[ y_m^- = \alpha(t_0-\tau_m) = \alpha(t_0) - \tau_m v^- + o(\tau_m). \]
Thus,
\[ d_F(y_m^+,y_m^-)
  \leq \lambda^+ |y_m^+-y_m^-|
  = \lambda^+\tau_m|v^++v^- + o(1) |. \]
For~$m$ large enough, we can take~$\lambda^+$ arbitrarily close to~$1$.
It follows from the inequality~$|v^+ + v^-| < 2$ that 
\[  d_F(y_m^-,y_m^+) < 2\tau_m  \]
contradicting that~$\alpha$ is a distance-realizing curve.
\end{proof}

Before proceeding to the proof of Theorem~\ref{thm:shortest_are_C1}, we extend the Finsler metric~$F$ to a surface~$M^+ \supseteq M$ with empty boundary; see Remark~\ref{rmk:extend}
As for any Finsler surface with empty boundary, every point of~$M^+$ has a normal neighborhood, that is, an open neighborhood~$U$ such that for any two points~$x,y\in U$, there is a unique geodesic from~$x$ to~$y$ contained in~$U$ and this geodesic is the unique distance-realizing arc from~$x$ to~$y$ in~$M^+$; see~\cite[p.~160]{bao2000introduction}. 
Note that if this geodesic is contained in~$M$, then it is also the unique distance-realizing arc from~$x$ to~$y$ in~$M$.

\begin{proof}[Proof of Theorem~\ref{thm:shortest_are_C1}]
We assume first that the metric is self-reverse.

Let~$\alpha:I \to M$ be a distance-realizing curve. Let~$t_0\in I$ and let~$x_0=\alpha(t_0)$.
If~$x_0=\alpha(t_0)$ lies in the interior of~$M$ then the arc~$\alpha$ coincides with a geodesic in a neighborhood of~$t_0$, where it is~$C^1$ (and we are done).
Thus, we can assume that~$x_0$ lies in~$\partial M$.

Again, we assume without loss of generality by working in a small enough neighborhood of~$x_0$ that~$M$ is a closed half-space of~$M^+=\R^2$ and that every geodesic arc is a unique distance-realizing arc.

Suppose that the arc~$\alpha$ is not differentiable on the right at some~$t_0\in I$.
(The argument is similar if~$\alpha$ is not differentiable on the left at~$t_0$.)
The arc~$\alpha$ has two departure velocities~$v$ and~$w$. Let~$K_v$ and~$K_w$ be two convex cones based at~$x_0$ that contain the points~$x_0+v$ and~$x_0+w$ in their interior and only meet at~$x_0$.
Take a unit vector~$u \in T_{x_0} M$ not tangent to~$\partial M$ that points in the interior of~$M$ and separates~$K_v$ from~$K_w$, and 
denote by~$\gamma_u$ the geodesic with initial velocity~$\gamma_u'(t_0)=u$. This geodesic does not visit~$K_v$ nor~$K_w$ in some interval~$(t_0,t_2)$.
On the other hand, the arc~$\alpha(t)$ visits the cones~$K_v$ and~$K_w$ infinitely many times in any interval $(t_0,\tau)$, with $\tau > t_0$. 
Therefore, it must cross the geodesic~$\gamma_u$ at some time~$t_1 \in (t_0,t_2)$. 
Since~$\gamma_u$ is the unique distance-realizing path between any of its points, the arc~$\alpha$ coincides with~$\gamma_u$ in~$[t_0,t_1]$. Thus~$\alpha$ does not visit~$K_v$ and~$K_w$ in~$(t_0,t_1)$. This contradiction proves that~$\alpha$ is differentiable on the right at~$t_0$. It follows from Lemma~\ref{lem:right} that~$\alpha$ is differentiable at every interior point~$t_0 \in I$.

Suppose~$\alpha$ is not~$C^1$ on the right at~$t_0$. (The argument is similar in case it is not~$C^1$ on the left.) 
The vector~$v=\alpha'(t_0)$ points inside~$M$ or is tangent to the boundary of~$M$. 
Since the velocities~$\alpha'(t)$ are unit vectors and the curve~$\alpha$ is not~$C^1$ on the right at~$t_0$, its derivative~$\alpha'$ has an accumulation point~$w\neq v$ when~$t$ goes to~$t_0$ from the right.
Let~$u$ be a unit vector spanning a line that separates~$v$ from~$w$.
Consider three disjoint neighbourhoods $U,\,V,\,W$ of~$u,\,v,\,w$ such that for every~$u',\,v',\,w'$ in~$U,\,V,\,W$ respectively, the line spanned by~$u'$ separates~$v'$ from~$w'$. Let~$K_V$ be the union of the rays contained in~$M$ starting at~$x_0$ with direction~$v' \in V$, and let~$R$ be any of these rays. Note that~$u$ is transverse to all these rays.
Working in a small enough neighbourhood of~$x_0$, we can assume that the family~$\Gamma$ of geodesics that visit~$R$ with velocity~$u$ foliates the cone~$K_V$, and that their tangent vectors do not deviate too much from~$u$ and thus lie in~$U$.
Since the velocity of~$\alpha$ at~$t_0$ lies in the open set~$V$, the arc~$\alpha$ restricted to some nontrivial interval $[t_0,t_3)$ lies in~$K_V$.
Now, since~$w$ is an accumulation point for~$\alpha'$ when~$t$ goes to~$t_0$ from the right, there exists $t_2 \in (t_0,t_3)$ such that $w' =\alpha'(t_2)$ lies in~$W$.
Let $x_2=\alpha(t_2)$, and let~$v'$ be the direction from~$x_0$ to~$x_2$. Let~$\gamma$ be the geodesic of~$\Gamma$ passing through~$x_2$, and let~$u'\in U$ be its velocity at~$x_2$.
The vector~$w'=\alpha'(t_2)$ points strictly inside the region of~$M$ delimited by~$\gamma$ containing~$x_0$, since the vector~$v'$ points outside, and the line generated by the vector~$u'$ separates~$v'$ from~$w'$.

Therefore, the arc~$\alpha$ starting at~$x_0$ must cross~$\gamma$ a first time at $t_1 \in (t_0,t_2)$ before crossing it again at~$t_2$.
Since~$\gamma$ is the unique distance-realizing path between $\alpha(t_1)$ and~$\alpha(t_2)$, the arc~$\alpha$ coincides with~$\gamma$ in~$[t_1,t_2]$, which contradicts the fact that~$\alpha'$ is transverse to~$\gamma$ at~$t_2$ (or~$t_1$). This finishes the proof of~\ref{thm:shortest_are_C1} for self-reverse metrics. 

In the case of directed metrics we adapt the argument as follows. Apart from the foliation~$\Gamma$, we need a second foliation~$\Gamma^-$ of~$K_V$ by geodesics transverse to the ray~$R$ with initial velocity~$-u$. Then we proceed as in the proof and after choosing the point~$x_1$ in~$K_V$, we let~$\gamma$ and~$\gamma^-$ be the two geodesics of~$\Gamma$ and~$\Gamma^-$ passing through~$x_1$. We keep only the part of each geodesic before it reaches~$x_1$ and discard the rest. These two half geodesics delimit a region of~$K_V$ containing~$x_0$.
The curve~$\alpha$ points strictly inside this region at~$x_1$.
Therefore, it must cross either~$\gamma$ or~$\gamma^-$ a first time before reaching~$x_1$.
We derive a contradiction as in the previous proof.
\end{proof}

\printbibliography
\end{document}